\documentclass[a4paper,12pt]{amsart}
\usepackage{amsmath}
\usepackage{amsfonts}
\usepackage{amssymb}
\usepackage{amsthm}
\usepackage[all]{xy}
\usepackage{geometry}
\usepackage{color}

\theoremstyle{definition}
\newtheorem{definition}{Definition}[section]
\newtheorem{remark}[definition]{Remark}

\newtheorem{example}[definition]{Example}

\theoremstyle{plain}
\newtheorem{theorem}[definition]{Theorem}

\newtheorem{lemma}[definition]{Lemma}
\newtheorem{proposition}[definition]{Proposition}
\newtheorem{corollary}[definition]{Corollary}

\DeclareMathOperator{\cA}{\mathcal{A}}
\DeclareMathOperator{\cB}{\mathcal{B}}
\DeclareMathOperator{\cC}{\mathcal{C}}
\DeclareMathOperator{\cD}{\mathcal{D}}

\DeclareMathOperator{\cM}{\mathcal{M}}
\DeclareMathOperator{\cN}{\mathcal{N}}
\DeclareMathOperator{\cO}{\mathcal{O}}
\DeclareMathOperator{\cK}{\mathcal{K}}

\DeclareMathOperator{\cS}{\mathcal{S}}
\DeclareMathOperator{\cT}{\mathcal{T}}

\DeclareMathOperator{\cX}{\mathcal{X}}
\DeclareMathOperator{\cY}{\mathcal{Y}}

\DeclareMathOperator{\add}{add}

\DeclareMathOperator{\Coker}{Coker}

\DeclareMathOperator{\coh}{coh}
\DeclareMathOperator{\cone}{cone}

\DeclareMathOperator{\depth}{depth}

\DeclareMathOperator{\End}{End}
\DeclareMathOperator{\gEnd}{\underline{End}}
\DeclareMathOperator{\Ext}{Ext}
\DeclareMathOperator{\gExt}{\underline{Ext}}

\DeclareMathOperator{\gr}{gr}

\DeclareMathOperator{\gldim}{gldim}

\DeclareMathOperator{\Gr}{Gr}

\DeclareMathOperator{\Hom}{Hom}
\DeclareMathOperator{\gHom}{\underline{Hom}}

\DeclareMathOperator{\id}{id}
\DeclareMathOperator{\idim}{idim}

\DeclareMathOperator{\im}{Im}
\DeclareMathOperator{\Ker}{Ker}

\DeclareMathOperator{\Lim}{Lim}

\DeclareMathOperator{\MCM}{MCM}

\DeclareMathOperator{\pdim}{pdim}

\DeclareMathOperator{\proj}{proj}
\DeclareMathOperator{\Proj}{Proj}
\DeclareMathOperator{\Qgr}{QGr}
\DeclareMathOperator{\qgr}{qgr}

\DeclareMathOperator{\tails}{tails}
\DeclareMathOperator{\Tor}{Tor}
\DeclareMathOperator{\tor}{tor}

\DeclareMathOperator{\D}{\mathbf{D}}

\begin{document}
\title{Noncommutative Resolutions of AS-Gorenstein isolated singularities}
\author{Haonan Li}
\address{School of Mathematics\\
Shanghai University of Finance and Ecnomics\\
 Shanghai, 200433\\
China}
\author{Menda Shen}
\author{Quanshui Wu}
\address{School of Mathematical Sciences \\
Fudan University \\
Shanghai, 200433 \\
 China}
\email{lihaonan@mail.shufe.edu.cn, 17110840001@fudan.edu.cn, qswu@fudan.edu.cn}
\thanks{This research has been supported by the NSFC (Grant No. 12471032) and the National Key Research and Development Program of China (Grant No. 2020YFA0713200).}

\keywords{Noncommutative resolution, noncommutative isolated singularity, noncommutative projective scheme, noncommutative projective coordinate ring, cluster tilting module}
\subjclass[2020]{14A22, 16S38, 16W50, 16E65, 14E15}

\newgeometry{left=3.18cm,right=3.18cm,top=2.54cm,bottom=2.54cm}

\begin{abstract}
In this paper, we investigate noncommutative resolutions of AS-Gorenstein isolated singularities. Noncommutative resolutions in graded case are achieved as the graded endomorphism rings of some finitely generated graded modules, which are seldom $\mathbb{N}$-graded algebras but bounded-below
$\mathbb{Z}$-graded algebras. So, the paper works on locally finite bounded-below
$\mathbb{Z}$-graded algebras.
We first define and study noncommutative projective schemes after Artin-Zhang, and define noncommutative quasi-projective spaces as the base spaces of noncommutative projective schemes. The equivalences between noncommutative quasi-projective spaces are proved to be induced by so-called modulo-torsion-invertible bimodules, which is in fact a Morita-like theory at the quotient category level. Based on the equivalences, we propose a definition of noncommutative resolutions of AS-Gorenstein isolated singularities, and prove that such noncommutative resolutions are generalized AS regular algebras.
The center of any noncommutative resolution is isomorphic to the center of the original AS-Gorenstein isolated singularity. In the final part we prove that a noncommutative resolution of an AS-Gorenstein isolated singularity of dimension $d$ is given by an MCM generator $M$ if and only if $M$ is a  $(d-1)$-cluster tilting module. A noncommutative version of the Bondal-Orlov conjecture is also proved to be true in dimension 2 and 3.
\end{abstract}

\maketitle

\tableofcontents

\section{Introduction}
In this paper we study resolutions of noncommutative isolated singularities in  noncommutative projective geometry.
Resolutions of singularities are important in algebraic geometry as they allow people to transfer calculations and constructions from singular varieties to smooth varieties. In algebraic geometry, a projective morphism $f: Y \to X$ is called a crepant resolution of $X$ if $Y$ is smooth and $f^*(\omega_X) = \omega_Y$ where $\omega_X$ and $\omega_Y$ are canonical bundles of $X$ and $Y$ respectively. Among resolutions of singularities, crepant resolutions share a good property as they preserve canonical modules. Bondal-Orlov \cite{BO1,BO2} conjectured that if $Y_1$ and $Y_2$ are crepant resolutions of a singular algebraic variety $X$, then the derived categories $\D^b(\coh Y_1)$ and $\D^b(\coh Y_2)$ are equivalent.
The conjecture is true in dimension 2 case \cite{KV}, and in some dimension 3 cases \cite{BKR, Bri}.
The bounded derived category of a crepant resolution of a singularity may be equivalent to the bounded derived category of some noncommutative algebra by tilting theory (see the example given in the introduction of \cite{V2}). McKay correspondence also gives a rich class of such noncommutative algebras. Based on this, Van den Bergh proposed a concept of noncommutative crepant resolutions (NCCRs for short) for a normal Gorenstein domain $R$ \cite{V1}. An NCCR of $R$ is an $R$-algebra $\Lambda$ of the form $\End_R(M)$ for some reflexive $R$-module $M$ such that $\Lambda$ has finite global dimension and $\Lambda$ is a maximal Cohen-Macaulay $R$-module.
Van den Bergh conjectured that all the (commutative and noncommutative) crepant resolutions of $R$ are derived equivalent. 
NCCRs have led to a lot of research both in algebraic geometry and noncommutative algebraic geometry.
The concept of NCCR was extended to some module-finite algebras, that is, algebras  finitely generated as modules over some central subrings. 
The NCCRs of module-finite algebras over complete regular local rings, normal Gorenstein domains or Cohen-Macaulay rings were studied respectively by Iyama \cite{I2}, Iyama-Reiten \cite{IR} and Iyama-Wemyss \cite{IW1,IW2}.

There is a well-known equivalence between the category of modules  over a commutative ring
$R$ and the category of quasi-coherent sheaves over the affine scheme Spec $R$. In noncommutative geometry, categories play the roles of geometric objects.
The philosophy of such categorical replacement is supported by many notable results.
One of these results is Gabriel-Rosenberg Reconstruction Theorem \cite{Ro}, which says that any scheme can be reconstructed uniquely up to isomorphism from the category of quasi-coherent sheaves on the scheme. 
Hence, all the information of a scheme is contained in its category of quasi-coherent sheaves.
Another is Serre's theorem \cite{Se}, which says that if $A$ is a finitely generated commutative graded algebra generated in degree one over a field $k$ and $X=\Proj A$ is the projective scheme associated to $A$, then there is an equivalence of categories
$F: \coh(X)\cong \qgr A$,
such that $F\circ (-\otimes \cO_X(n))\cong s^n F$ and $F$ sends the structure sheaf $\cO_X$ to $\cA$, where $\qgr A$ is the quotient category of the finitely generated graded module category modulo finite dimensional modules, $\cA$ is the image of $A$ in $\qgr A$ (for notations, see \S \ref{Quotient Category}).

An algebraic triple $(\cC,\cA,s)$ is considered as a noncommutative projective scheme by Artin-Zhang \cite{AZ} where $\cC$ is an abelian category, $\cA$ is an object in $\cC$ and $s$ is an auto-equivalent functor of $\cC$. Here $\cA$ is regarded as the noncommutative structure sheaf, and $s$ as the polarization of the noncommutative projective scheme $(\cC,\cA,s)$. In particular, in \cite{AZ}, $(\qgr A,\cA,s)$ is called the noncommutative projective scheme associated to $A$ for a noetherian $\mathbb{N}$-graded $k$-algebra $A$, and $A$ is regarded as the coordinate ring of this noncommutative projective scheme. Artin-Zhang proved a noncommutative version of Serre's theorem, which gives a necessary and sufficient condition for an algebraic triple $(\cC,\cA,s)$ being equivalent to the associated noncommutative projective scheme of a noetherian $\mathbb{N}$-graded ring $A$ (see Theorem \ref{noncommutative Serre theorem}).

A modified definition of noncommutative projective schemes is given in Definition \ref{def-noncom-proj-scheme} as an algebraic triple satisfying some additional conditions. 
It is worth pointing out that the definition of a noncommutative projective scheme contains the specific noncommutative structure sheaf and polarization (an auto-equivalent functor).
This differs from the commutative case. The Gabriel-Rosenberg Reconstruction Theorem suggests that taking the category of quasi-coherent sheaves will not lose any information of the scheme, so there is no necessary to emphasize the structure sheaf once the category is given. In fact, an abelian category may admit different noncommutative structure sheaves and/or different polarizations giving different noncommutative projective schemes, see Example \ref{projective space P1}.

As an analogue of isolated singularities in commutative case, noncommutative isolated singularities are studied in \cite{Jo, Ue1, Ue2, MU,LW2} for $\mathbb{N}$-graded algebras.
The work of NCCRs aforementioned provides a kind of noncommutative resolutions of some isolated singularities. Inspired by their work, we study the noncommutative resolutions of noncommutative AS-Gorenstein isolated singularities.
Note that in projective geometry, the projective scheme associated to an isolated singularity is smooth, so it has no meaning to consider this question in commutative realm (see also Theorem \ref{main theorem: Morita} (5)).

A resolution is to find some smooth object which is equivalent to a singular one in some sense.
In \cite{QWZ}, Qin-Wang-Zhang used the equivalence of quotient categories of module categories with respect to some dimension function to define noncommutative quasi-resolutions for Auslander-Gorenstein algebras and such noncommutative quasi-resolutions are Auslander-regular algebras.
Enlightened by \cite{QWZ}, He-Ye \cite{HY} investigated the resolution of noncommutative isolated singularities. They defined a right quasi-resolution of noncommutative isolated singularity $A$ to be a graded endomorphism ring $B$ of a maximal Cohen-Macaulay $A$-module $M$ such that $B$ is a (generalized) AS-regular algebra and $\gHom_A(M,-)$ induces an equivalence from $\qgr A$ to $\qgr B$.

AS-regular algebras are introduced by Artin and Schelter in 1987 \cite{AS}, as an analogue of polynomial rings in noncommutative setting.  AS-regular algebras
are regarded as the coordinate rings of noncommutative projective spaces. 
Thus the noncommutative projective schemes associated to AS-regular algebras are the smooth objects in our mind.

As in \cite{QWZ} and \cite{HY}, that a resolution is birational equivalent to the singularity could be replaced by the equivalence of some quotient categories.
Therefore, roughly speaking, a noncommutative resolution of an AS-Gorenstein isolated singularity $A$ is a (generalized) AS-regular algebra $B$ such that there is an equivalent functor
$F:\qgr A \cong \qgr B$ and $Fs\cong sF$
(written as $F:(\qgr A,s)\cong (\qgr B,s)$).
Given a noncommutative projective scheme $(\qgr A,\cA,s)$, $(\qgr A,s)$ is called the noncommutative quasi-projective space of $(\qgr A,\cA,s)$ (see Definition \ref{def-noncom-proj-scheme}). If $A$ is a (generalized) AS-regular algebra, then $(\qgr A,s)$ is called a noncommutative projective space.

(Noncommutative) McKay correspondence provides a key inspiration for the development of noncommutative resolutions of AS-Gorenstein isolated singularities. For a (skew) polynomial ring $A$ and a group action $G$ on $A$, under some conditions, the skew group algebra $A\#G$ is generalized AS-regular and $(\qgr A^G,s)\cong (\qgr A\#G, s)$. 
These facts lead us to consider the equivalent functor $F:\qgr A\cong \qgr B$  commuting with $s$ or say $F$ preserves the noncommutative polarization.
There is a very nice review paper \cite{Leu} by Leuschke about the relations between noncommutative resolutions and McKay correspondence in the commutative case. For the noncommutative case, one can see \cite{CYZ, MU} and Example \ref{example of Zhu and CYZ}. 
We mention that it is also valuable to consider the equivalent functor $F':\qgr A\cong \qgr B$ preserving the noncommutative structure sheaf, which we will not consider in the present paper. There are some related works in \cite{ATV,AV}.

Noncommutative resolutions are often achieved as the endomorphism ring of some modules. For noncommutative isolated singularities, their noncommutative resolutions are graded endomorphism rings of some finitely generated graded modules. Such kind of graded endomorphism rings may not be $\mathbb{N}$-graded, but being bounded-below $\mathbb{Z}$-graded.
Locally finite bounded-below $\mathbb{Z}$-graded algebras are called commonly graded in \cite{CKWZ2}.
This suggests us that to study resolutions of noncommutative isolated singularities we should work on commonly graded algebras instead of $\mathbb{N}$-graded algebras. Most basic results for $\mathbb{N}$-graded algebras \cite{LW,RR} are in fact true for commonly graded algebras.

In this paper we mainly focus on noncommutative resolutions of commonly graded AS-Gorenstein algebras (see Definition \ref{generalized-AS-Gorenstein}, sometimes they are called generalized AS-Gorenstein algebras in literature)  which are noncommutative isolated singularities (see Definition \ref {def-nc-iso-sing}). Noncommutative Serre's theorem \cite[Theorem 4.5, Corollary 4.6]{AZ} plays a key role in the study of noncommutative algebraic geometry. As a preparation, we modify \cite[Theorem 4.5, Corollary 4.6]{AZ} to cover the commonly graded case. We follow the notations in \cite{AZ}.

\begin{theorem}[Corollary \ref{coor. ring}]
\begin{enumerate}
        \item [(1)] Suppose $(\cC,O,s)$ is a noncommutative projective scheme (see Definition \ref{def-noncom-proj-scheme}). Then $B=B(\cC,O,s)$ is a right noetherian commonly graded algebra satisfying $\chi_1$, $\depth_B(B)\geqslant 2$, and the algebraic triples
$(\cC,O,s)$ and $(\qgr B,\cB,s)$ are isomorphic.
        \item [(2)] If $A$ is a right noetherian commonly graded algebra such that $\depth_A(A)\geqslant 2$, then $A$ satisfies $\chi_1$ if and only if that $(\qgr A, \cA, s)$ is a noncommutative projective scheme. In this case, $A$ is isomorphic to $B(\qgr A, \cA, s)$.
    \end{enumerate}
\end{theorem}

To find smooth noncommutative projective schemes resolving noncommutative AS-Gorenstein isolated singularities, a crucial step is to characterize the equivalences between noncommutative quasi-projective spaces. 
Such equivalences, 
similarly to the Morita equivalence of module categories,  are given by a kind of graded bimodules, which are called modulo-torsion-invertible bimodules (see Definition \ref{definition of modulo-torsion-invertible}).

\begin{theorem}[Theorem \ref{morita theory}, Corollary \ref{center of noncommutative coordinate ring}]\label{main theorem: Morita}
Let $A$ and $B$ be two noncommutative projective coordinate rings which are noetherian.
Suppose there is an equivalence of noncommutative quasi-projective spaces
\[F:(\qgr A,s)\rightleftarrows (\qgr B,s): G\]
Let $M=\omega G\cB$ and $M'=\omega F\cA$. Then
\begin{enumerate}
    \item [(1)] $B\cong \gEnd_A(M)$ and $A\cong \gEnd_B(M')$.
    \item [(2)] ${}_BM_A$ and ${}_AM'_B$ are finitely generated on both sides.
    \item [(3)] $F\cong \pi_B\gHom_{\cA}(\cM,-)\cong -\otimes_{\cA}\cM'$ and $G\cong\pi_A\gHom_{\cB}(\cM',-)\cong -\otimes_{\cB}\cM$.
    \item [(4)] $M\cong\gHom_B(M',B)$ and $M'\cong \gHom_A(M,A)$ as graded bimodules.\\ Thus there is a graded Morita context $(A,B,M,M',\tau,\mu)$ isomorphic to the graded Morita context defined by $M_A$ or $M'_B$.
    \item [(5)] The center of $A$ and the center of $B$ are isomorphic.
\end{enumerate}
\end{theorem}

The module $M$ in Theorem \ref{main theorem: Morita} is exactly a modulo-torsion-invertible $(B,A)$-bimodule $M$ with $\depth_A(M)\geqslant 2$  (see Theorem \ref{modulo-torsion-invertible bimodule and equivalent functor}).
Moreover, the equivalences of noncommutative quasi-projective spaces $(\qgr A,s)$ and $(\qgr B,s)$ are in one-to-one correspondence to the modulo-torsion-invertible $(B,A)$-bimodules with depth no less than $2$ up to isomorphism (see Theorem \ref{modulo-torsion-invertible bimodule and equivalent functor}).
Theorem \ref{main theorem: Morita} also shows that the equivalence is induced by the graded Morita context defined by $M_A$ and $B$ is isomorphic to $\gEnd_A(M)$.
Therefore, a noncommutative resolution of a commonly graded AS-Gorenstein isolated singularity should be a graded endomorphism ring of some modulo-torsion-invertible module with depth no less than $2$.

Let $A$ be a commonly graded AS-Gorenstein isolated singularity.
By the Morita-like theory of noncommutative quasi-projective spaces, we investigate the maximal Cohen-Macaulay (MCM for short) generators $M_A$ over $A$ and its graded endomorphism ring $B=\gEnd_A(M)$. In fact, such an $M$ is a modulo-torsion-invertible $(B,A)$-bimodule. Based on this, we define a noncommutative resolution of $A$ to be a graded endomorphism ring $B$ of an MCM generator $M_A$ such that $B$ has finite global dimension and $B_B$ is MCM (see Definition \ref{definition of NCCR}). In fact, $B$ is a commonly graded AS-regular algebra such that the noncommutative projective space $(\qgr B,s)$ associated to $B$ is equivalent to the noncommutative quasi-projective space $(\qgr A,s)$ associated to $A$.

\begin{theorem}[Proposition \ref{tails of NCCR equi. to tails of singularity}, Theorem \ref{NCCR is AS regular}]
    Let $A$ be a commonly graded AS-Gorenstein isolated singularity of dimension $d\geqslant 2$.  
    Then the following are equivalent.
    \begin{enumerate}
        \item [(1)] $B$ is a noncommutative resolution of $A$ given by $M_A$.
        \item [(2)] $B^o$ is a noncommutative resolution of $A^o$ given by ${}_AM'=\gHom_A(M,A)$ and $M_A$ is an MCM module.
        \item [(3)] $B=\gEnd_A(M)$ is a noetherian commonly graded AS-regular algebra of dimension $d$ and $M_A$ is an MCM generator.
    \end{enumerate}
If (1), (2) or (3) holds, then the graded Morita context defined by $M_A$ induces an equivalence:
\[-\otimes_{\cA}\cM':(\qgr A,s)\rightleftarrows (\qgr B,s):-\otimes_{\cB}\cM.\]
\end{theorem}

Similar to the case of module-finite algebras \cite{I2,IR}, every noncommutative resolution of a commonly graded AS-Gorenstein isolated singularity $A$ of dimension $d$ is given by a $(d-1)$-cluster tilting $A$-module (see Definition \ref{cluster tilting module}); and conversely, every $(d-1)$-cluster tilting $A$-module gives  a noncommutative resolution of $A$.

\begin{theorem}[Theorem \ref{d-1 CT module and NCCR}]
Suppose $A$ is a commonly graded AS-Gorenstein isolated singularity of dimension $d\geqslant 2$.  Then the following are equivalent.
\begin{enumerate}
    \item [(1)] $M$ is a $(d-1)$-cluster tilting $A$-module.
    \item [(2)] $B=\gEnd_A(M)$ is  a noncommutative resolution of $A$ given by $M_A$.
\end{enumerate}
\end{theorem}

In the last part of the paper, we prove that a version of noncommutative Bondal-Orlov conjecture holds for the noncommutative resolutions of commonly graded AS-Gorenstein isolated singularities in dimension $2$ and $3$.

\begin{theorem}[Theorem \ref{proof for BO conjecture}]
Suppose $A$ is a commonly graded AS-Gorenstein isolated singularity of dimension $d$.
\begin{enumerate}
    \item [(1)] If $d=2$, then all noncommutative resolutions of $A$ are Morita equivalent.
    \item [(2)] If $d=3$, then all noncommutative resolutions of $A$ are derived Morita equivalent.
\end{enumerate}
\end{theorem}

The paper is organized as follows.
In section \ref{Preliminaries}, we introduce some notations and some related basic material. In section \ref{Noncommutative Projective Schemes}, we  define noncommutative projective schemes and noncommutative quasi-projective spaces following \cite{AZ}, and study their coordinate rings. In section \ref{Morita Theory for Noncommutative quasi-projective spaces}, a Morita-like theory of noncommutative quasi-projective spaces is given. It is proved that the equivalences of noncommutative quasi-projective spaces are in one-to-one correspondence to the modulo-torsion-invertible bimodules with depth no less than $2$.
In section \ref{NC-singu}, we study the commonly graded AS-Gorenstein isolated singularities, and prove that all its MCM generators are modulo-torsion-invertible.
In section \ref{Noncommutative Crepant Resolutions}, we show that all the noncommutative resolutions are commonly graded AS-regular algebras. We also prove that a noncommutative resolution of a AS-Gorenstein isolated singularity of dimension $d$ is given by an MCM generator $M$ if and only if $M$ is a  $(d-1)$-cluster tilting module. In dimension $2$ and $3$, a corresponding version of the noncommutative Bondal-Orlov conjecture is proved at the end.

\section{Preliminaries}\label{Preliminaries}
\subsection{Notations and Conventions}
Let $A$ be a $\mathbb{Z}$-graded $k$-algebra over a fixed field $k$.
A graded $A$-module will always mean a right graded $A$-module. Let $A^o$ be the opposite algebra of $A$ and $A^e=A\otimes_k A^o$. Then a graded $A^o$-module is exactly a left graded $A$-module and a graded $A^e$-module is exactly a graded $(A,A)$-bimodule.

Let $M(n)$ be the $n$-th shift of a graded module $M$ where $M(n)_i =M_{n+i}$. Let $X[n]$ be the $n$-th shift of a complex $X$ with $X[n]^i=X^{n+i}$.
The category of graded $A$-modules is denoted by $\Gr A$ and the full subcategory of finitely generated graded $A$-modules is denoted by $\gr A$. The derived category of $\Gr A$ is denoted by $\D(\Gr A)$.

Let $D(-)=\gHom_k(-,k)$ be the Matlis duality.
A noetherian algebra $A$  means that it is left and right noetherian.

For a graded $A$-module $M$, let $\add_A M$ be the full subcategory of $\Gr A$  consisting of all objects $X$ such that $X$ is isomorphic to a direct summand of a finite direct sum of the shifts of $M$.

\subsection{Commonly Graded Algebras}
As mentioned, a noncommutative resolution should be the graded endomorphism ring of some finitely generated graded module. Such a graded endomorphism ring is usually not $\mathbb{N}$-graded but commonly graded. 
Some basic results about commonly graded algebras are collected here for convenience (see
\cite{NO} and \cite{CKWZ2}).

\begin{definition} 
    A $\mathbb{Z}$-graded $k$-algebra $A = \oplus_{i \in \mathbb{Z}} A_i$  is called \textit{locally finite} if $A_i$ is finite-dimensional for all $i\in \mathbb{Z}$. A locally finite $\mathbb{Z}$-graded algebra is called \textit{commonly graded} if $A$ is bounded-below, that is, $A_{\leqslant n} := \oplus_{i \leqslant n} A_i = 0$ for some $n \in \mathbb{Z}$.
\end{definition}

Most of the results for $\mathbb{N}$-graded algebras hold for commonly graded algebras. We will cite some related results directly if their proofs in $\mathbb{N}$-graded case still work for commonly graded algebras  without any essential change. If $M$ is a graded module over a commonly graded algebra $A$, $M_{\geqslant n}:= \oplus_{i \geqslant n} M_i$ may not be a graded $A$-submodule. But under some finiteness conditions, as showed in Lemmas \ref{J and A>i cofinal} and \ref{J^nM and M>n cofinal}, the submodule $M_{\geqslant n}$ used in the proofs for $\mathbb{N}$-graded algebras can be replaced by the submodule $MJ_A^n$,
where $J_A$ is the graded Jacobson radical of $A$.

For a commonly graded algebra $A$, the \textit{graded Jacobson radical} of $A$, denoted by $J_A$ or simply $J$, is the intersection of all maximal graded left ideals. Note that $J$ is also the intersection of all maximal graded right ideals of $A$ \cite[Prop. 2.9.1]{NO}.

Let
$b_l(M)=\inf\{i\mid M_i\neq 0\}\text{ and }b_u(M)=\sup\{i\mid M_i\neq 0\}$ for a  $\mathbb{Z}$-graded $k$-space $M \neq 0$.

\begin{lemma}\label{A-J finite-dim}
\cite[Lemma 3.2]{CKWZ2} Let $A$ be a commonly graded algebra.
\begin{enumerate}
\item[(1)] If $L$ is a maximal graded left ideal of $A$, then $L\supseteq A_{> -b_l(A)}$. As a consequence, $J\supseteq A_{>-b_l(A)}$.
\item[(2)] For any maximal graded left ideal $L$  of $A$,
\[b_l(A/L)\geqslant b_l(A) \textrm{ and } b_u(A/L)\leqslant -b_l(A).\]
As a consequence, $A/L$ is finite-dimensional.
\item[(3)] $A/J$ is finite-dimensional and is graded semisimple, i.e. $A/J$ is a direct sum of graded simple modules.
\end{enumerate}
\end{lemma}

\begin{lemma}\label{J^i is contained in A>j} For any $i>0$, there is some integer $n_i$ such that $J^{n_i}\subseteq A_{\geqslant i}$. As a consequence, $\bigcap_nJ^n=0$.
\end{lemma}
\begin{proof}
Let $I=AA_{\geqslant 1-2b_l(A)}A$, which is an ideal of $A$ contained in $A_{\geqslant 1}$. Let $B=A/I$. Then $B$ is finite-dimensional. Hence $J_B^k=0$ for some integer $k$. Thus $J^k\subseteq I\subseteq A_{\geqslant 1}$. Therefore $J^{ik}\subseteq A_{\geqslant i}$ for all $i$. 
\end{proof}

The following is a graded version of Nakayama's Lemma \cite[Proposition 3.3(2)]{CKWZ2} \cite[Lemma 1.7.5]{NO}.
\begin{lemma}\label{Nakayama's Lemma}
If $M$ is a bounded-below graded $A$-module such that $MJ=M$, then $M=0$.
\end{lemma}
\begin{proof}
For a fixed $i > 0$, by Lemma \ref{J^i is contained in A>j}, there is some $n$ such that  $J^n\subseteq A_{\geqslant i}$. Then $M=MJ=MJ^n\subseteq MA_{\geqslant i}$. If $M \neq 0$, it follows from the fact that $M$ is bounded-below that  $b_l(M)\geqslant b_l(M) + i$, which is a contradiction.
\end{proof}

The following lemma is the commonly graded version of \cite[Lemma 2.3]{RR}.

\begin{lemma}\label{finitely generated algebra}
 $J/J^2$ is finite-dimensional if and only if ${}_AS=A/J$ (or equivalently, $S_A$) is finitely presented. In this case, $A$ is finitely generated as $k$-algebra.
\end{lemma}
\begin{proof}
By Lemma \ref{A-J finite-dim}, $A/J$ is finite-dimensional. Hence, $J/J^2$ is finite-dimensional if and only if  $J$ is finitely generated as a left (right) $A$-module by Lemma \ref{Nakayama's Lemma}.
This is equivalent to that ${}_AS$ (or equivalently, $S_A$) is a finitely presented module.

Let $V$ be a finite-dimensional subspace of $J$ such that $V+J^2=J$, and $X$ be a finite-dimensional space of $A$ such that $X+J=A$. By induction, $V^n+J^{n+1}=J^n$. Then $A=X+V+V^2+\cdots+V^n+J^{n+1}$ for all $n$. By lemma \ref{J^i is contained in A>j}, for any $i > 0$ there is an integer $n_i$ such that $J^{n_i}\subseteq A_{\geqslant i}$. Thus $A_{i-1}\cap J^{n_i}=0$. This implies $A_{i-1}\subseteq X+V+V^2+\cdots+V^{n_i-1}$. Therefore $A$ is generated as a $k$-algebra by $A_{\leqslant 0}$, $V$ and $X$.
\end{proof}

\begin{lemma}\label{J and A>i cofinal}
If $J/J^2$ is finite-dimensional, then $\{J^n\}_n$ and $\{A_{\geqslant n}\}_n$ are cofinal.
\end{lemma}
\begin{proof}
By Lemma \ref{J^i is contained in A>j}, we only need to prove that for any $i>0$ there is an integer $j$ such that $A_{\geqslant j}\subseteq J^i$.

With the notation as in the proof of Lemma \ref{finitely generated algebra},
$A=X+V+V^2+\cdots+V^{i-1}+J^i.$
Since $V$ and $X$ are finite-dimensional, $X+V+V^2+\cdots+V^{i-1}$ are both bounded-below and bounded-above. Thus for $j\gg 0$, $A_{\geqslant j}\subseteq J^i$.
\end{proof}

\begin{lemma}\label{J^nM and M>n cofinal}
If $J/J^2$ is finite-dimensional, then for any finitely generated graded $A$-module $M$, $\{M_{\geqslant n}\}_n$ and $\{MJ^n\}_n$ are cofinal.
\end{lemma}
\begin{proof}
By Lemma \ref{J and A>i cofinal}, it suffices to show that $\{M_{\geqslant n}\}_n$  and  $\{MA_{\geqslant n}\}_n$ are cofinal.

Since $M$ is finitely generated, $M$ is bounded-below. Hence, for any $n$, there is an integer $m$ such that $MA_{\geqslant m} \subseteq M_{\geqslant n}$.

On the other hand, let $\{m_1,\cdots, m_s\}$ be a set of homogeneous generators of $M$ and  $r=\max\{r_1,\cdots, r_s\}$ where $r_i=\deg m_i$. Then, for any $n$ and  any homogeneous element  $m\in M$ with $\deg m \geqslant n+r$, $m=\Sigma m_ia_i$ for some $a_i \in A$ with $\deg a_i=\deg m-r_i$. So,  $m\in MA_{\geqslant n}$. Hence, $M_{\geqslant n+r}\subseteq MA_{\geqslant n}$.
\end{proof}

\subsection{Local cohomology for commonly graded algebras}
Let $A$ be a commonly graded algebra, and $J$ be its graded Jacobson radical.
For a graded $A$-module $M$, let
\[\Gamma_A(M):=\{m\in M\mid m\cdot J^i =0, i\gg 0\}.\]
Then $\Gamma_A(M)\cong  \underrightarrow{\Lim}\gHom_A(A/J^i,M)$.
 If $J/J^2$ is finite-dimensional, then by Lemma \ref{J and A>i cofinal},
$\Gamma_A(M)=\{m\in M\mid m\cdot A_{\geqslant i}=0, i\gg 0\}$. If $\Gamma_A(M)=M$ (resp. $\Gamma_A(M)=0$), then $M$ is called \textit{torsion} (resp. \textit{torsion-free}).

The \textit{$n$-th local cohomology} of $M$ is
\[R^n\Gamma_A(M)\cong \underrightarrow{\Lim}\gExt_A^n(A/J^i,M).\]
The \textit{cohomological dimension} of $\Gamma_A$ is
\[\sup\{i\mid R^i\Gamma_A(M)\neq 0,M\in \Gr A\}.\]

Next proposition is well-known for connected graded algebras \cite[Thm 8.3]{AZ}.

\begin{proposition}\label{local cohomology over B and over A when B_A finite}
Let $A$ and $B$ be commonly graded algebras, and $f:A\to B$ be a graded algebra morphism such that ${}_AB$ is finitely generated. Suppose $A$ is noetherian and $B$ is right noetherian. Then for any $N\in \D^+(\Gr B)$ there is a canonical isomorphism
\[R\Gamma_A(N) \cong R\Gamma_B(N).\]
Furthermore, if both $\Gamma_A$ and $\Gamma_B$ have finite cohomological dimension, then the same result holds for any $N\in \D(\Gr B)$.
\end{proposition}
\begin{proof} We first claim that $\Gamma_A(Y)=\Gamma_B(Y)$ for any graded $B$-module $Y$.
Since ${}_AB$ is finitely generated, $B/(J^n_A\cdot B)\cong (A/J^n_A)\otimes_AB $ is finite-dimensional. Thus $J^t_B\subseteq J^n_A\cdot B$ for some $t$.
If  $y\in Y$ such that $0 =y\cdot J_A^n \in Y_A$, then $yJ_B^t\subseteq y\cdot(J^n_A\cdot B)=0$.
On the other hand,
if $yJ_B^n=0$, then by Lemma \ref{J and A>i cofinal} there exist integers $n_1$ and $n_2$ such that
\[y\cdot J_A^{n_2}\subseteq y\cdot A_{\geqslant n_1}=yf(A_{\geqslant n_1})\subseteq yB_{\geqslant n_1}\subseteq yJ_B^n=0.\]
Therefore $y \in \Gamma_A(Y)$ if and only if $y \in \Gamma_B(Y)$. Hence
$\Gamma_A(Y)=\Gamma_B(Y).$

Next we claim that any graded injective $B$-module $I$ is $\Gamma_A$-acyclic, that is, $R^i\Gamma_A(I)=0$ for all $i>0$. Since $B$ is right noetherian, we may assume $I_B$ is either torsion or torsion-free by \cite[Lemma 5.1]{LW}.

If $I_B$ is torsion then $I_A$ is torsion. Hence  there is a graded injective resolution $E^{\bullet}$ of $I_A$ such that all $E^i_A$ are torsion. So, $R^i\Gamma_A(I)=\textrm{H}^i(E^{\bullet})=0$ for all $i>0$.

If $I_B$ is torsion-free, then $\gHom_B(X,I)=0$ for any finite-dimensional module $X_B$.
Let $(P_{\bullet})_A$ be a graded projective resolution of $A/J^n_A$. Then
\[\xymatrix{
 \gHom_B(P_{i-1}\otimes_A B,I)\ar[r]\ar[d]^{\cong}& \gHom_B(P_{i}\otimes_A B,I)\ar[r]\ar[d]^{\cong}&\gHom_B(P_{i+1}\otimes_AB ,I)\ar[d]^{\cong} \\
\gHom_A(P_{i-1},I)\ar[r] &\gHom_A(P_{i},I)\ar[r] &\gHom_A(P_{i+1},I).
}\]
It follows that
$\gHom_B(\Tor^A_{i}(A/J_A^n,B),I)\cong \gExt_A^i(A/J_A^n,I)$.

Since ${}_AB$ is finitely generated, ${}_AB$ has a finitely generated free resolution which implies that $\Tor^A_{i}(A/J_A^n,B)$ is a finite-dimensional $B$-module. So $\gHom_B(\Tor^A_{i}(A/J_A^n,B),I)=0$
 and $\gExt_A^i(A/J_A^n,I)=0$ for all $i$. Hence $R\Gamma_A(I)=0$.

Now, for any $N\in \D^+(\Gr B)$, and a graded injective $B$-module resolution $I^{\bullet}$ of $N$, $R\Gamma_B(N)=\Gamma_B(I^{\bullet})=\Gamma_A(I^{\bullet})=R\Gamma_A(N)$.

Suppose that $\Gamma_A$ and $\Gamma_B$ have finite cohomological dimension. For any $N \in \D(\Gr B)$, $N$ has a $\Gamma_B$-acyclic resolution, which is also a $\Gamma_A$-acyclic resolution of $N$ by the above argument.  It follows that $R\Gamma_B(N)\cong R\Gamma_A(N)$.
\end{proof}

For any graded $A$-module $M$,  the \textit{depth} of $M$ is defined by \[
\depth_A(M)=\inf\{i\mid \gExt_A^i(A/J,M)\neq 0\}.
\]

The proofs of the following lemmas are the same as in the $\mathbb{N}$-graded case.

\begin{lemma}\label{depth and local cohomology}
    For a graded $A$-module $M$,
    \[\inf\{i\mid \gExt_A^i(A/J,M)\neq 0\}=\inf\{i\mid R^i\Gamma_A(M)\neq 0\}.\]
\end{lemma}

\begin{lemma}\label{depth in exact sequences}
    For any exact sequence
    $0\to M\to X_{s-1}\to X_{s-2}\to \cdots\to X_0\to N\to 0$ of graded $A$-modules,
    \[\depth_A(M)\geqslant \inf\{\depth_A(N)+s,\depth_A(X_0)+s-1,\cdots,\depth_A(X_{s-1})\}.\]
\end{lemma}

\subsection{Quotient Category}\label{Quotient Category} We mainly refer to \cite[Chapter 4]{Po}
 for the general quotient category theory.
Let $\cC$ be an abelian category and $\cS$ a Serre subcategory of $\cC$. The quotient category of $\cC$ with respect to $\cS$ is denoted by $\cC/\cS$. Let $\pi: \cC\to \cC/\cS$ be the quotient functor.

In general, $\cC/\cS$ is an abelian category and $\pi$ is exact. If $\pi$ has a right adjoint functor $\omega$, then $\cS$ is called a \textit{localizing subcategory} of $\cC$. Note that in this case $\pi\omega\cong \id_{\cC/\cS}$ and $\omega$ is fully faithful.

Next proposition has been proved in \cite[Corollary 15.9]{Fa}, but the assumptions in the corollary are not enough. It is fixed in \cite[Proposition 2.1.4]{L} and \cite[Lemma 1.1]{S}.

\begin{proposition}\label{induced functors on quotient categories}
    Let  $F: \cC \to \cD$ be a functor between two abelian categories, $\cS$ resp. $\cT$ be a Serre subcategory of $\cC$ resp. $\cD$.
    \begin{enumerate}
        \item [(1)] Suppose that $\cC$ has enough projective objects and $F$ is right exact. Then there is a  functor $F':\cC/\cS \to \cD/\cT$ such that $\pi F\cong F'\pi$ if and only if $F(\cS)\subseteq \cT$ and $L^1F(\cS)\subseteq \cT$. The functor $F'$ is unique up to natural isomorphism.  If $\cS$ is a localizing subcategory of $\cC$, then $F'\cong \pi F\omega$.
        \item [(2)] Suppose that  $\cC$ has enough injective objects and $F$ is left exact. Then there is a  functor $F':\cC/\cS \to \cD/\cT$ such that $\pi F \cong F'\pi$ if and only if $F(\cS)\subseteq \cT$ and $R^1F(\cS)\subseteq \cT$. The functor $F'$ is unique up to natural isomorphism. If $\cS$ is a localizing subcategory of $\cC$, then $F'\cong \pi F\omega$.
    \end{enumerate}
\[
\xymatrix{
\cC \ar[r]^F \ar@<-1mm>[d]_{\pi} & \cD \ar[d]^{\pi} \\
\cC/\cS \ar@{.>}[r]^{F'} \ar@<-1mm>[u]_{\omega} & \cD/\cT
}\]
When $\pi F\cong F'\pi$, we say $F$ induces $F'$ or $F'$ is induced by $F$.
\end{proposition}

Let $A$ and $B$ be two right noetherian commonly graded algebras and ${}_BM_A$ be a graded $(B,A)$-bimodule. Suppose $\cS$ and $\cT$ are Serre subcategories of $\gr A$ and $\gr B$ respectively. Then $M$ is called \textit{compatible} with $\cS$ and $\cT$ if $-\otimes_{B}M$ induces a functor
\[-\otimes_{\cB}\cM:\gr B/\cT\to \gr A/\cS.\]
By Proposition \ref{induced functors on quotient categories}, $M$ is compatible with $\cS$ and $\cT$ if and only if for any $N\in\cT$, both $N\otimes_{B}M$ and $\Tor^{B}_1(N,M)$ are in $\cS$.

 We  refer to \cite{Ja} for the notations related to  Morita theory. For a graded Morita context $(A,B,M,M',\tau,\mu)$, if $\tau$ is surjective, then $\tau$ is an isomorphism. Next lemma is a generalization of this fact.

\begin{lemma}\cite[Lemma 2.2]{L}\label{coker tau in cS}
Suppose $(A,B,M,M',\tau,\mu)$ is a graded Morita context where $A$ is right noetherian, $M_A$ and $M'_{B}$ are finitely generated. If $\cS$ is a Serre subcategory of $\gr A$, then $\Coker\tau\in \cS$ implies $\Ker \tau \in \cS$.
\end{lemma}
\begin{proof}
    Since $\Ker\tau\cdot \im \tau=0$, $\Ker\tau$ can be regarded as a finitely generated right $(A/\im\tau)$-module. Hence $\Ker\tau\in \cS$.
\end{proof}

The following theorem shows when a Morita context induces an equivalence between the quotient categories.

\begin{theorem}\cite[Theorem 2.2.5]{L}\label{equ. of quot. cat. induced by Morita context}
Suppose $(A,B,M,M',\tau,\mu)$ is a graded Morita context where $A$ and $B$ are right noetherian, $M_A$ and $M'_{B}$ are finitely generated. Let $\cS$ and $\cT$ be Serre subcategories of $\gr A$ and $\gr B$ respectively such that $M$ and $M'$ are compatible with them. Then
\[
-\otimes_{\cA}\cM':\gr A/\cS \rightleftarrows \gr B/\cT :-\otimes_{\cB}\cM
\]
is an equivalence if and only if $\Coker\tau\in \cS$ and $\Coker \mu\in \cT$.
\end{theorem}
\begin{proof}
Suppose $\Coker\tau\in \cS$ and $\Coker \mu\in \cT$. By Lemma \ref{coker tau in cS}, $\Ker \tau\in \cS$ and $\Ker \mu\in \cT$.
Consider the following two exact sequences
    \[0\to \Ker \tau\to M'\otimes_{B}M\to \im \tau\to 0,\]
    \[0\to \im \tau\to A\to \Coker\tau\to 0.\]

For any $N\in \gr A$, $N\otimes_A-$ induces the following exact sequences
    \[\cdots\to N\otimes_A\Ker \tau\to N\otimes_A(M'\otimes_{B}M)\to N\otimes_A\im \tau\to 0,\]
    \[\cdots\to \Tor_1^A(N,\Coker \tau)\to N\otimes_A\im \tau\to N\otimes_A A\to N\otimes_A\Coker\tau\to 0.
    \]

Since $\Ker\tau\in\cS$ and $\cS$ is a Serre subcategory of $\gr A$, 
 $N\otimes_A\Ker \tau, N\otimes_A\Coker\tau$ and $\Tor_1^A(N,\Coker \tau)$ are in $\cS$. Hence in $\gr A/\cS$,
\[\cN\otimes_{\cA}\cM'\otimes_{\cB}\cM \cong \pi(N\otimes_A(M'\otimes_{B}M))\cong \pi(N\otimes_A\im\tau)\cong \pi(N\otimes_AA)\cong \cN.\]

Similarly, for any $N'\in\gr B$,
\[\cN'\otimes_{\cB}\cM\otimes_{\cA}\cM' \cong \pi(N'\otimes_{B}(M\otimes_AM'))\cong \cN'.\]
So $
-\otimes_{\cA}\cM':\gr A/\cS \rightleftarrows \gr B/\cT :-\otimes_{\cB}\cM
$ is an equivalence between the  quotient categories.

The inverse direction holds obviously.
\end{proof}

When the conditions in Theorem \ref{equ. of quot. cat. induced by Morita context} hold, an equivalence between the quotient categories
\[F:\gr A/\cS \rightleftarrows \gr B/\cT :G\]
is called \textit{induced by the Morita context} $(A,B,M,M',\tau,\mu)$ if $F\cong -\otimes_{\cA}\cM'$ and $G\cong -\otimes_{\cB}\cM$.

\begin{remark}
    If $A,B$ are not necessarily right noetherian and $M_A,M'_{B}$ are not necessarily finitely generated, we consider $\Gr A$ and $\Gr B$ instead of  $\gr A$ and $\gr B$, and assume $\cS$ and $\cT$ are  localizing subcategories. Then an analogous result of Theorem \ref{equ. of quot. cat. induced by Morita context} still holds. Ungraded versions of previous results also hold.
\end{remark}

Now we consider the Serre category generated by torsion modules. Let $\Tor A$ be the full subcategory consisting of graded torsion $A$-modules and $\tor A=\Tor A\cap \gr A$. Let
\[\Qgr A=\Gr A/\Tor A \text { and } \qgr A=\gr A/\tor A.\]
The quotient functor is denoted by $\pi_A$ (or simply $\pi$). Let $\omega_A$ (or simply $\omega$) be the right adjoint of $\pi:\Gr A\to \Qgr A$. Let $s$ denote the auto-equivalent functor of $\Qgr A$ induced by shift functor $(1)$ of $\Gr A$.

For a graded $A$-module $M$, $\pi M$ is denoted by $\cM$.
For any $\cM,\cN\in \Qgr A$, let
\[\Hom_{\cA}(\cM,\cN)=\Hom_{\Qgr A}(\cM,\cN) \text{ and } \gHom_{\cA}(\cM,\cN)=\bigoplus_{i\in\mathbb{Z}}\Hom_{\cA}(\cM,s^i\cN).\]

Since $\Qgr A$ has enough injective objects which are the images of torsion-free injective graded $A$-modules, the derived functor of $\Hom_{\cA}(\cM,-)$ exists.
Let
\[\Ext_{\cA}^i(\cM,\cN)=\Ext_{\Qgr A}^i(\cM,\cN) \text{ and } \gExt_{\cA}^i(\cM,\cN)=\bigoplus_{i\in\mathbb{Z}}\Ext_{\cA}^i(\cM,s^i\cN).\]

Next, let $\cS=\tor A$ and $\cT=\tor B$.

\begin{lemma}\label{induced functor on tails by tensor}
Let $(A,B,M,M',\tau,\mu)$ be a graded Morita context where $A,B$ are right noetherian commonly graded algebras and $M_A,M'_B$ are finitely generated. Then $M$ (resp. $M'$) is compatible with $\tor A$ and $\tor B$ if and only if ${}_BM$ (resp. ${}_AM'$) is finitely presented.
\end{lemma}
\begin{proof}
Suppose $M$ is compatible with $\tor A$ and $\tor B$.
Let $P_1\to P_0\to {}_BM\to 0$ be a minimal projective resolution of ${}_BM$. Then $B/J_B\otimes_BM\cong B/J_B\otimes_BP_0$ and $\Tor^B_1(B/J_B,M)\cong B/J_B\otimes_BP_1$ are finite-dimensional. Thus both $P_0$ and $P_1$ are finitely generated and ${}_BM$ is finitely presented.

Suppose ${}_BM$ is finitely presented, and $0\to K\to P\to {}_BM\to 0$ is an exact sequence of finitely generated graded left $B$-modules where ${}_BP$ is free. For any $Y\in \tor B$, it induces a long exact sequence
\[0\to \Tor^B_1(Y,M)\to Y\otimes_BK\to Y\otimes_BP\to Y\otimes_BM\to 0.\]
Then $Y\otimes_BP$ is finite-dimensional. Hence $Y\otimes_BM$ is finite-dimensional. Similarly, $Y\otimes_BK$ is finite-dimensional, so $\Tor_1^B(Y,M)$ is finite-dimensional. Therefore $M$ is compatible with $\tor A$ and $\tor B$.
\end{proof}

Here is the version of Theorem \ref{equ. of quot. cat. induced by Morita context} for $\cS=\tor A$ and $\cT=\tor B$.

\begin{theorem}\label{equ. of quot. cat. induced by Morita context for commonly graded}
Let $(A,B,M,M',\tau,\mu)$ be a graded Morita context where $A,B$ are right noetherian commonly graded algebras and $M_A,M'_B$ are finitely generated. Suppose ${}_BM$ and ${}_AM'$ are finitely presented. Then the Morita context $(A,B,M,M',\tau,\mu)$ induces an equivalence:
\[-\otimes_{\cA}\cM':\qgr A \rightleftarrows \qgr B :-\otimes_{\cB}\cM \]
 if and only if $\Coker \tau$ and $\Coker\mu$ are finite-dimensional.
\end{theorem}

For convenience, we introduce a notion of modulo-torsion-invertible $(B,A)$-bimodule.

\begin{definition}\label{definition of modulo-torsion-invertible}
Let $A$ and $B$ be noetherian commonly graded algebras. If there is a graded Morita context $(A,B,{}_BM_A,{}_AM'_B,\tau,\mu)$ with ${}_BM_A$ and ${}_AM'_B$ finitely generated on both sides such that $\Coker \tau$ and $\Coker \mu$ are finite-dimensional, then ${}_BM_A$ is called a \textit{modulo-torsion-invertible $(B,A)$-bimodule} with the inverse ${}_AM'_B$ associated to $(A,B,{}_BM_A,{}_AM'_B,\tau,\mu)$.
\end{definition}

\subsection{commonly graded AS-Gorenstein algebras and MCM modules}
Commonly graded ($\mathbb{N}$-graded) AS-Gorenstein (AS-regular) property is studied in \cite{LW} and \cite{RR}. Most properties still hold in the commonly graded case. We recall the definitions first.

\begin{definition}\label{generalized-AS-Gorenstein}
A commonly graded algebra $A$ is called \textit{Artin-Schelter Gorenstein} (for short, AS-Gorenstein) of dimension $d$ if the following conditions hold.
\begin{enumerate}
\item[(1)] $A$ has  finite left and right graded injective dimension $d$.
\item[(2)] For every graded simple $A$-module $M$, $\gExt^i_A(M,A)=0$ if $i\neq d$; for every graded simple $A^o$-module $N$, $\gExt_{A^o}^i(N,A)=0$ if $i\neq d$.
\item[(3)] $\gExt_A^d(-,A)$ and $\gExt_{A^o}^d(-,A)$ give a bijection between the isomorphic classes of graded simple $A$-modules to the isomorphic classes of graded simple $A^o$-modules.
\end{enumerate}
If moreover $A$ has finite graded global dimension $d$, then $A$ is called \textit{Artin-Schelter regular} (for short, AS-regular).

If $A$ is $\mathbb{N}$-graded, then it is called an \textit{$\mathbb{N}$-graded AS-Gorenstein (regular) algebra}.
\end{definition}


As we extend our sight from $\mathbb{N}$-graded algebras to commonly graded algebras, it is natural to ask the existence of commonly graded AS-Gorenstein (regular) algebras with non-zero negative part. The answer is positive.
\begin{example}
    Let $A$ be an $\mathbb{N}$-graded AS-Gorenstein (regular) algebra and $B=\gEnd_A(A\oplus A(1))$. Then $B$ is a commonly graded algebra with $B_{<0}$ non-zero.
    By \cite[Proposition 6.2]{LW}, the property of $\mathbb{N}$-graded AS-Gorenstein is preserved by graded Morita equivalence for $\mathbb{N}$-graded algebras. This is also true for commonly graded algebras with a same proof.
    It follows that $B$ is a commonly graded AS-Gorenstein (regular) algebra.
    \end{example}

An important property of noetherian commonly graded AS-Gorenstein algebras is the existence of balanced dualizing complexes. Balanced dualizing complex is defined as follows.

\begin{definition}
    Suppose $A$ is a noetherian commonly graded algebra. A complex $R\in \D^b(\Gr A^e)$ is called a \textit{dualizing complex} of $A$, if it satisfies the following conditions:
\begin{enumerate}
\item[(1)] $R$ has finite injective dimension over $A$ and $A^o$ respectively.
\item[(2)] The homologies of $R$ are finitely generated as $A$-modules and $A^o$-modules.
\item[(3)] The natural morphisms $ A\to R\gHom_A(R,R)$ and $A\to \gHom_{A^o}(R,R)$ are isomorphisms in $\D(\Gr A^e)$.
\end{enumerate}
If moreover, $R\Gamma_A(R)\cong D(A)$ and $R\Gamma_{A^o}(R)\cong D(A)$ in $\D(\Gr A^e)$, then $R$ is called a \textit{balanced dualizing complex} of $A$.
\end{definition}

\begin{theorem}\label{AS-G has balanced dualizing complex and chi}
    Let $A$ be a noetherian commonly graded AS-Gorenstein algebra of dimension $d$.
    \begin{enumerate}
        \item [(1)] $A$ admits a balanced dualizing complex given by $D(R^d\Gamma_A(A))[d]$.
        \item [(2)] $A$ (resp. $A^o$) satisfies the condition  $\chi$ (see Definition \ref{def-chi}).
        \item [(3)] $\Gamma_A$ (resp. $\Gamma_{A^o}$) has finite cohomological dimension $d$.
    \end{enumerate}
\end{theorem}
\begin{proof}
    By \cite[Theorem 5.12, Corollary 5.13]{LW}.
\end{proof}

Next we define maximal Cohen-Macaulay modules.

\begin{definition}
    Let $A$ be a commonly graded algebra. Suppose $\Gamma_A$ has finite cohomological dimension $d$. A finitely generated graded $A$-module $M$ is called a \textit{maximal Cohen-Macaulay} (MCM for short) $A$-module if $\depth_A (M)=d$.
\end{definition}

The full subcategory consisting of all MCM $A$-modules is denoted by $\MCM (A)$.

\begin{lemma}\label{depth < infinity}
    Let $A$ be a noetherian commonly graded algebra with a balanced dualizing complex.
    Then $\depth_A(M)<\infty$ for any $0\neq M\in\Gr A$.
\end{lemma}
\begin{proof}
By \cite[Lemma 3.6]{LW}.
\end{proof}

\begin{definition}
    A noetherian commonly graded algebra $A$ is called a \textit{balanced Cohen-Macaulay algebra} of dimension $d$ if it admits a balanced dualizing complex $R$ such that $R\cong D(R^d\Gamma_A(A))[d]$ in $\D(\Gr A^e)$.
\end{definition}

We introduce some well-known results about MCM modules in the following.

\begin{proposition}\label{Hom(-,A) induce duality on MCM}
    Let $A$ be a noetherian commonly graded AS-Gorenstein algebra of dimension $d$ and $M$ be a finitely generated graded $A$-module.
    \begin{enumerate}
        \item [(1)] $M$ is MCM if and only if $\gExt_A^i(M,A)=0$ for all $i>0$.
        \item [(2)] There is a duality
    \[\gHom_A(-,A):\MCM (A)\rightleftarrows \MCM (A^o):\gHom_{A^o}(-,A).\]
    \end{enumerate}
\end{proposition}
\begin{proof}
    (1) If $A$ is commonly graded AS-Gorenstein, then $D(R^d\Gamma_A(A))$ is an invertible graded $(A,A)$-bimodule by a commonly graded version of \cite[Proposition 5.8]{LW}. 
    It follows that $\gExt_A^i(M,A)=0$ if and only if $\gExt_A^i(M,D(R^d\Gamma_A(A)))=0$. By Local duality, 
    \[R^i\Gamma_A(M)\cong D(\gExt_A^{d-i}(M,D(R^d\Gamma_A(A))).\]
    Hence $M$ is MCM if and only if $\gExt_A^i(M,A)=0$ for all $i>0$ and $\gHom_A(M,A)\neq 0$.

    Suppose $\gExt_A^i(M,A)=0$ for all $i>0$. Since $A$ has finite injective dimension on both sides, the double Ext spectral sequence converges
    \[E^{p,q}_2=\gExt_{A^o}^p(\gExt_A^{-q}(M,A),A)\Rightarrow\left\{
    \begin{aligned}
        0,\quad &p+q\neq 0,\\
        M,\quad & p+q=0.
    \end{aligned}
    \right.
    \]
    It follows that $\gHom_{A^o}(\gHom_A(M,A),A)\cong M$. Thus $\gHom_A(M,A)\neq 0$. Therefore $M$ is MCM if and only if $\gExt_A^i(M,A)=0$ for all $i>0$.

    (2) It follows from (1).
\end{proof}

\begin{proposition}\label{MCM and syzygy}
    Let $A$ be a noetherian commonly graded AS-Gorenstein algebra and $M$ an MCM $A$-module. Then for any $N\in \gr A$ and $i,j>0$,
    \[\gExt_A^i(M,N)\cong \gExt_A^{i+j}(M,\Omega^jN)\]
    where $\Omega$ is the syzygy functor.
\end{proposition}
\begin{proof}
    Let $P_\bullet\to N$ be a finitely generated projective resolution of graded $A$-module $N$. Then for any $j$, there is an exact sequence,
    \[0\to \Omega^{j+1}N\to P_j\to \Omega^jN\to 0.\]
    It induces a long exact sequence
    \begin{align*}
        0&\to \gHom_A(M,\Omega^{j+1}N)\to \gHom_A(M,P_j)\to \gHom_A(M,\Omega^jN)\\
        &\to \gExt_A^1(M,\Omega^{j+1}N)\to\gExt_A^1(M,P_j)\to\gExt_A^1(M,\Omega^jN)\to\cdots.
    \end{align*}

    Since $M$ is MCM, $\gExt_A^i(M,P_j)=0$ for all $i>0$ by Proposition \ref{Hom(-,A) induce duality on MCM}. So
    \[\gExt_A^i(M,\Omega^jN)\cong \gExt_A^{i+1}(M,\Omega^{j+1}N)\]
    for all $i>0$.  Then the assertion follows from an induction.
\end{proof}
An important class of MCM modules is cluster tilting modules. The concept of cluster tilting is introduced in higher Auslander-Reiten theory \cite{I1,I2} (it is called maximal orthogonal originally) and it can be regarded as a generalization of Cohen-Macaulay-finiteness, that is, the isomorphic classes of MCM $A$-modules are finite up to shift.

\begin{definition} \label{cluster tilting module}
    Let $A$ be a commonly graded algebra such that $\Gamma_A$ has finite cohomological dimension. An MCM $A$-module $M$ is called an \textit{$n$-cluster tilting module} if
    \begin{align*}
        \add_A M&=\{N\in \MCM (A)\mid \gExt_A^i(M,N)=0,\forall\, 0<i<n\}\\
        &=\{N\in \MCM (A)\mid \gExt_A^i(N,M)=0,\forall\, 0<i<n\}.
    \end{align*}
\end{definition}

By definition, $A$ is Cohen-Macaulay-finite if and only if $A$ has a $1$-cluster tilting module.

\section{Noncommutative projective schemes}\label{Noncommutative Projective Schemes}
In this section, we define and study noncommutative projective schemes and their coordinate rings after \cite{AZ}. An \textit{algebraic triple} $(\cC,O,s)$, which consists of a $k$-linear abelian category $\cC$, an object $O\in \cC,$ and a $k$-linear auto-equivalence $s$ of $\cC$, is called a noncommutative projective scheme if it satisfies the conditions $H1$, $H2$, $H3$ in \cite{AZ}, and an additional condition $H4$ (see Theorems \ref{noncommutative Serre theorem} and \ref{noncommutative Serre theorem for commonly graded algebra}). We prove that $(\cC,O,s)$ is a noncommutative projective scheme if and only if the algebraic triple $(\cC,O,s)$ is isomorphic to $(\qgr A,\cA,s)$ for some right noetherian commonly graded algebra $A$ satisfying $\chi_1$ and $\depth_A(A)\geqslant 2$ (see Corollary \ref{coor. ring}).

\subsection{Noncommutative projective schemes after Artin-Zhang}
The condition $\chi$ plays a very important role in noncommutative projective geometry, which is introduced in \cite{AZ} for $\mathbb{N}$-graded algebras. The $\chi$-condition can be defined similarly for commonly graded algebras.

\begin{definition}\label{def-chi}
Let $A$ be a commonly graded $k$-algebra, $M$ be a graded $A$-module.
\begin{enumerate}
    \item [(1)] If $\dim_k\gExt^j_A(A/J,M)<\infty$  for all $j\leqslant i$, then we say that \textit{$\chi_i(M)$ holds}.

    \item [(2)] If $\chi_i(M)$ holds for every finitely generated graded $A$-module $M$, then we say that \textit{$\chi_i$ holds for $A$} or \textit{$A$ satisfies the  $\chi_i$-condition}.

    \item [(3)] If  $\chi_i$ holds for $A$ for all $i \geqslant 0$, then we say  \textit{$\chi$ holds for $A$} or \textit{$A$ satisfies the  $\chi$-condition}.
\end{enumerate}
\end{definition}

The results in \cite[Section 3]{AZ} hold for right noetherian commonly graded algebras under some modifications, say, replacing $A_{\geqslant n}$ by $J^n$, etc. We collect some results here for later use.

\begin{proposition}\label{facts about chi condition}
    Let $A$ be a right noetherian commonly graded algebra.
\begin{enumerate}
    \item [(1)] $\chi_i$ holds for $A$ if and only if $R^j\Gamma_A(M)$ is bounded-above for all $j\leqslant i$ and all $M\in \gr A$.
    \item [(2)] $\chi_1$ holds for $A$ if and only if $(\omega\!\cM)_{\geqslant n}$ is contained in some finitely generated graded submodule of $\omega\!\cM$ for any $M\in \gr A$ and $n\in\mathbb{Z}$.
    \item [(3)] If $A$ satisfies $\chi_i$, then for any finitely generated graded $A$-modules $M$, $N$ and $j< i$, the canonical map from $\gExt_A^j(M,N)$ to $\gExt_{\cA}^j(\cM,\cN)$ has bounded-above kernel and cokernel.
\end{enumerate}
\end{proposition}
\begin{proof}
    The proofs are similar to that of \cite[Corollary 3.6(3), Proposition 3.14 and  Corollary 7.3]{AZ} in the $\mathbb{N}$-graded case respectively.
\end{proof}

 If $A$ is right noetherian commonly graded, then $\Qgr A$ and $\qgr A$ have similar cohomological theory as the  $\mathbb{N}$-graded case  given in \cite[Section 7]{AZ}.

\begin{proposition}\label{cohomology in tails}
    Suppose that $A$ is a right noetherian commonly graded algebra and $M$ is a graded $A$-module. Then
    \begin{enumerate}
    \item [(1)]
    $R^{i+1}\Gamma_A(M)\cong \gExt_{\cA}^i(\cA,\cM)$ for any $i\geqslant 1$.
    \item [(2)] there is a long exact sequence
    \[0\to \Gamma_A(M)\to M\xrightarrow[]{\varphi_M} \omega\pi M\to R^1\Gamma_A(M)\to 0.\]
    \end{enumerate}
\end{proposition}
\begin{proof}
Similar to the proof of \cite[Proposition 7.2]{AZ} in the $\mathbb{N}$-graded case.
\end{proof}

\begin{corollary}\label{depth geqslant 2}
 \begin{enumerate}
    \item [(1)] $\depth_A(M)\geqslant 2$ if and only if $\varphi_M:M\to \omega\pi M$ is an isomorphism.
    \item [(2)] $\depth_A(\omega\pi M)\geqslant 2$.
    \item [(3)] If $\depth_A(M)\geqslant 2$, then $\gHom_A(M,M)\cong \gHom_{\cA}(\cM,\cM)$ as $\mathbb{Z}$-graded algebras.
     \item [(4)] If $\depth_A(A)\geqslant 2$, then $A\cong\gHom_{\cA}(\cA,\cA)$ as commonly graded algebras.
 \end{enumerate}
\end{corollary}
\begin{proof}
    (1) By Lemma \ref{depth and local cohomology}, $\depth_A(M)\geqslant 2$ if and only if $\Gamma_A(M)=R^1\Gamma_A(M)=0$.

(2) and (3) follow from (1). 

(4) follows from (3).
\end{proof}

\begin{definition}\cite{AZ}\label{conditions of ample} Let $(\cC,O,s)$ be an algebraic triple.
The pair $(O,s)$ is called \textit{ample} in $\cC$ if
\begin{enumerate}
    \item [(1)] for every object $\cM \in \cC$, there are some positive integers $r_1,r_2,\cdots,r_p$ and an epimorphism $\oplus_{i=1}^ps^{-r_i}O\to\cM$ in $\cC$;
    \item [(2)] for every epimorphism $\cM\to\cN$ in $\cC$, there is an integer $n_0$ such that the induced map $\Hom_{\cC}(s^{-n}O,\cM)\to\Hom_{\cC}(s^{-n}O,\cN)$ is surjective for every $n\geqslant n_0$.
\end{enumerate}
\end{definition}

Without loss of generality, the auto-equivalence $s$ of $\cC$ can be assumed to be an automorphism (see \cite[Proposition 4.2]{AZ}).

Let $(\cC,O,s)$ be an algebraic triple. Then
$B(\cC,O,s):=\oplus_{i\in \mathbb{Z}}\Hom_{\cC}(O,s^iO)$
is a $\mathbb{Z}$-graded algebra with a natural multiplication, and $\oplus_{i\in \mathbb{Z}}\Hom_{\cC}(O,s^i\cM)$ has a natural graded right $B(\cC,O,s)$-module structure for any $\cM\in\cC$.

Let $A$ be a (right noetherian) commonly graded algebra and $\cC=\Qgr A$ ($\cC = \qgr A$). Then for any $\cM\in \cC$, $(\cC,\cM,s)$ is an algebraic triple and $B(\cC,\cM,s)=\gHom_{\cA}(\cM,\cM)$ is a $\mathbb{Z}$-graded algebra, where $s$ is the shift functor.

Two algebraic triples $(\cC_1,O_1,s_1)$ and $(\cC_2,O_2,s_2)$ are called \textit{isomorphic} if there is an equivalent functor $F:\cC_1\to \cC_2$ such that $FO_1\cong O_2$ and $s_2F\cong Fs_1$. If $G:\cC_1\to \cC_2$ is an equivalent functor such that $s_2G\cong Gs_1$, then $(O_1,s_1)$ is ample in $\cC_1$ if and only if $(GO_1,s_2)$ is ample in $\cC_2$.

Next theorem is called noncommutative Serre's theorem in literature,  as it is a noncommutative analogue of Serre's theorem in algebraic geometry.

\begin{theorem}(See \cite[Theorem 4.5, Corollary 4.6]{AZ})\label{noncommutative Serre theorem}
\begin{enumerate}
\item [(1)] Let $(\cC,O,s)$ be an algebraic triple. Suppose
   \begin{enumerate}
     \item [H1.] $O\in \cC$ is a noetherian object;
     \item [H2.] $\Hom_{\cC}(O,\cM)$ is finite-dimensional for all $\cM\in \cC$;
     \item [H3.] $(O,s)$ is ample in $\cC$.
   \end{enumerate}

Then $B=B(\cC,O,s)_{\geqslant 0}$ is a right noetherian locally finite $\mathbb{N}$-graded  $k$-algebra satisfying $\chi_1$, and the functor
\[\cC\to \qgr B,\,\, \cM\mapsto \pi(\oplus_{i\in Z}\Hom_{\cC}(O,s^i\cM))\]
gives an isomorphism between  algebraic triples
$(\cC,O,s)$ and $ (\qgr B,\cB,s)$.
\item [(2)] Let $A$ be a right noetherian locally finite $\mathbb{N}$-graded  algebra satisfying $\chi_1$. Then H1, H2 and H3 in (1) hold for the algebraic triple $(\qgr A, \cA, s)$.
Moreover, there is a canonical morphism $A\to B$ for $B=B(\qgr A,\cA,s)_{\geqslant 0}$ with the kernel and cokernel being  bounded-above, so it induces an isomorphism  between $(\qgr A,\cA, s)$ and $(\qgr B,\cB,s)$.
\end{enumerate}
\end{theorem}

The following lemma is extracted and modified from the proof of \cite[Theorem 4.5]{AZ}) (see Theorem \ref{noncommutative Serre theorem}).
\begin{lemma}\label{facts extracted from noncom. Serre thm}
 Suppose the algebraic triple $(\cC,O,s)$ satisfies $H1$, $H2$, $H3$ in Theorem \ref{noncommutative Serre theorem}, and satisfies
   \begin{enumerate}
     \item [H4.] $B(\cC,O,s)$ is bounded-below.
   \end{enumerate}
    Let $B=B(\cC,O,s)$ and $\Phi=\oplus_{i\in Z}\Hom_{\cC}(O,s^i-):\cC\to \Gr B$. Then for any $\cM\in \cC$,
    \begin{enumerate}
        \item [(1)] $\Phi(\cM)_{\geqslant n}$ is contained in a finitely generated $B$-submodule of $\Phi(\cM)$ for any $n\in \mathbb{Z}$;
        \item [(2)] the morphism $\varphi$ in the following exact sequence
        \[0\to \Gamma_B(\Phi(\cM))\to \Phi(\cM)\xrightarrow[]{\varphi} \omega\pi \Phi(\cM)\to R^1\Gamma_B(\Phi(\cM))\to 0\]
        is an isomorphism.
    \end{enumerate}
\end{lemma}
\begin{proof}
    (1) Similar to the proof of S5 in \cite[Theorem 4.5]{AZ}.

    (2) Similar to the proof of the corollary of S7 in \cite[Theorem 4.5]{AZ}.
\end{proof}

\begin{corollary}\label{depth of coordinate ring of noncom. proj. scheme}
    Let $(\cC,O,s)$ be as in Lemma \ref{facts extracted from noncom. Serre thm}. Then, $B=B(\cC,O,s)$ satisfies $\chi_1$, and $\depth_B(\Phi(\cM))\geqslant 2$ for any $\cM\in \cC$.
    In particular, $\depth_B(B)\geqslant 2$.
\end{corollary}

The following theorem is modified for commonly graded algebras from Theorem \ref{noncommutative Serre theorem} with almost the same proof. It can be viewed as a generalized version of noncommutative Serre's theorem.
\begin{theorem}\label{noncommutative Serre theorem for commonly graded algebra}
\begin{enumerate}
\item [(1)] Let $(\cC,O,s)$ be an algebraic triple. Suppose it satisfies $H1$, $H2$, $H3$ in Theorem \ref{noncommutative Serre theorem}, and satisfies
   \begin{enumerate}
     \item [H4.] $B(\cC,O,s)$ is bounded-below.
   \end{enumerate}

Then $B:=B(\cC,O,s)$ is a right noetherian commonly graded $k$-algebra  satisfying $\chi_1$ and $\depth_B(B)\geqslant 2$. The functor
\[\cC\to \qgr B,\,\,\cM\mapsto \pi(\oplus_{i\in Z}\Hom_{\cC}(O,s^i\cM))\]
gives an isomorphism between algebraic triples
$(\cC,O,s)$ and $(\qgr B,\cB,s)$.
\item [(2)] Let $A$ be a right noetherian commonly graded algebra satisfying $\chi_1$. Then $H1$, $H2$, $H3$ hold for $(\qgr A, \cA, s)$, and there is a canonical morphism $A\to B$ such that both the kernel and cokernel of it are bounded-above, where $B=B(\qgr A,\cA,s)$.

    If $B=B(\qgr A,\cA,s)$ is bounded-below (this is true when $\depth_A(A)\geqslant 2$ by Corollary \ref{depth geqslant 2}), then there is an isomorphism  between $(\qgr A,\cA, s)$ and  $(\qgr B,\cB,s)$.
\end{enumerate}
\end{theorem}

Now we are ready to give a definition for noncommutative quasi-projective spaces and noncommutative projective schemes.

\begin{definition}\label{def-noncom-proj-scheme}
    Let $\cC$ be a $k$-linear abelian category and $s$ be a $k$-linear auto-equivalence of $\cC$. If, for some object $O$ in $\cC$, the algebraic triple $(\cC,O,s)$ satisfies $H1$, $H2$, $H3$ and $H4$ in Theorem \ref{noncommutative Serre theorem for commonly graded algebra}, then $(\cC,O,s)$ is called a \textit{noncommutative projective scheme}, $O$ is called its \textit{structure sheaf} and $s$ is called its \textit{polarization}. The pair $(\cC,s)$ is called the \textit{noncommutative quasi-projective space} of $(\cC,O,s)$.
\end{definition}

Noncommutative projective schemes defined in Definition \ref{def-noncom-proj-scheme} are stronger than the (noetherian) projective schemes defined in \cite{AZ}, where the hypothesises $H1$, $H2$, $H3$ and $H4$ are not assumed.

Two noncommutative projective schemes are called \textit{isomorphic}  if they are isomorphic  as algebraic triples.
Two noncommutative quasi-projective spaces $(\cC_1,s_1)$ and $(\cC_2,s_2)$ are called \textit{equivalent} if there is an equivalent functor $F:\cC_1\to \cC_2$ such that $s_2F\cong Fs_1$. It follows from the definition that for any structure sheaf $O_1$ of $(\cC_1,s_1)$, $FO_1$ is a structure sheaf of $(\cC_2,s_2)$ and $F$ induces an isomorphism between noncommutative projective schemes $(\cC_1,O_1,s_1)$ and $(\cC_2,FO_1,s_2)$.

On one hand, for any noncommutative projective scheme $(\cC,O,s)$, $B(\cC,O,s)$ is a right noetherian commonly graded algebra by Theorem \ref{noncommutative Serre theorem for commonly graded algebra}.
Moreover, $B(\cC,O,s)$ satisfies $\chi_1$ and its depth $\geqslant 2$ by Corollary \ref{depth of coordinate ring of noncom. proj. scheme}.  On the other hand, if $A$ is a right noetherian commonly graded algebra satisfying $\chi_1$ and $\depth_A(A)\geqslant 2$, then $(\qgr A, \cA, s)$ is a noncommutative projective scheme by Theorem \ref{noncommutative Serre theorem for commonly graded algebra} and  Corollary \ref{depth geqslant 2}.
Consequently, we have the  following corollary.

\begin{corollary}\label{coor. ring}
\begin{enumerate}
        \item [(1)] Suppose $(\cC,O,s)$ is a noncommutative projective scheme. Then $B=B(\cC,O,s)$ is a right noetherian commonly graded algebra satisfying $\chi_1$, $\depth_B(B)\geqslant 2$, and the algebraic triples
$(\cC,O,s)$ and $(\qgr B,\cB,s)$ are isomorphic.
        \item [(2)] If $A$ is a right noetherian commonly graded algebra such that $\depth_A(A)\geqslant 2$, then $A$ satisfies $\chi_1$ if and only if that $(\qgr A, \cA, s)$ is a noncommutative projective scheme. In this case, $A$ is isomorphic to $B(\qgr A, \cA, s)$.
    \end{enumerate}
\end{corollary}

\begin{remark}\label{abelian category admits different structure sheaves and quto-equi}
In general, there is no Gabriel-Rosenberg Reconstruction Theorem for noncommutative projective scheme. In fact, as Example \ref{projective space P1} shows,
    a noncommutative quasi-projective space $(\cC,s)$ may admit different structure sheaves.
    Given an abelian category $\cC$ and an object $O\in \cC$, there may exist different auto-equivalences $s_1$ and $s_2$ such that $(\cC,O,s_1)$ and $(\cC,O,s_2)$ are noncommutative projective schemes which are not isomorphic.
 \end{remark}

\begin{example}\label{projective space P1}
    Let $A=k[x_0,x_1]$ be the graded polynomial algebra with $\deg x_0=\deg x_1 = 1$, and $X=\Proj A$ be the projective scheme associated to $A$ with the structure sheaf $\cO_X$.
    By Serre's theorem,
    \[(\coh X,\cO_X,-\otimes \cO_X(1))\cong (\qgr A,\cA,s).\]

    Let $B$ be the $2$-Veronese subring of $A$. The subring $B$ is viewed as a graded algebra with $B_i=A_{2i}$. In fact, $B$ is isomorphic to $k[y_0,y_1,y_2]/(y_0y_2-y_1^2)$ as graded algebras, where $\deg  y_i=1$ for $i  = 0, 1, 2$. Let $Y=\Proj B$. By \cite[Proposition 5.10]{AZ} and Serre's theorem 
    \[(\coh X,\cO_X,-\otimes \cO_X(2))\cong (\qgr A,\cA,s^2)\cong (\qgr B,\cB,s)\cong(\coh Y,\cO_Y,-\otimes \cO_Y(1)).\]

    Let $\sigma$ be the automorphism of $A$ given by $\sigma(x_0)=-x_0,\sigma(x_1)=-x_1$, and  $G$ be the cyclic group generated by $\sigma$. 
    Then the invariant subring is $A^G=k[x_0^2,x_0x_1,x_1^2]$. Let $Z=\Proj A^G$. By the Remark after \cite[Proposition 5.10]{AZ}, 
    \[(\qgr A^G,\pi(A^G),s)\cong (\qgr B\oplus \qgr B,(\cB,0),\tilde{s})\]
    where $\tilde{s}$ is induced by an auto-equivalence of $\gr B\oplus \gr B$ which sends $(M,N)$ to $(N,M(1))$. Serre's theorem doesn't apply to $A^G$, as $A^G$ is not generated in degree one. 
    Although the gradations of the Veronese subring $B$ and the invariant subring $A^G$ are different, $Y = \Proj B \cong \Proj A^G = Z$.
    Then
    \[(\coh Z,\cO_Z,-\otimes \cO_Z(2))\cong(\coh Y,\cO_Y,-\otimes \cO_Y(1)) \cong (\qgr B,\cB,s).\]
    So, $(\coh Z,\cO_Z,-\otimes \cO_Z(2))$ is not isomorphic to $(\qgr A^G,\pi(A^G),s)$.

    Let $A\#G$ be the skew group algebra of the $G$-action. By \cite[Corollary 3.11]{MU}, there is an isomorphism
    \[(\qgr A\#G, \pi(A\#G),s)\cong (\qgr A^G,\cA,s)\]
    where $\cA$ is the image of $A_{A^G}$ in $\qgr A^G$. Then one sees  that $(\qgr A^G,\pi(A^G),s)$ and $(\qgr A^G,\cA,s)$ are not isomorphic. This offers an example that there are noncommutative projective schemes such that the abelian categories and the polarizations are the same, while the noncommutative structure sheaves are not isomorphic.

    In fact, as projective schemes, $X$, $Y$ and $Z$ are exactly the projective space $\mathbb{P}_k^1$. But as noncommutative projective schemes, $(\coh X,\cO_X,-\otimes \cO_X(1))$ and $(\coh Y,\cO_Y,-\otimes \cO_Y(1))$ are different. Note that $(\coh Y,\cO_Y,-\otimes \cO_Y(1))$ is isomorphic to $ (\coh X,\cO_X,-\otimes \cO_X(2))$. This offers an example that there are different noncommutative projective schemes with the same abelian categories and the same structure sheaves, but the polarizations are not naturally isomorphic.
\end{example}

\begin{remark}
    Let $A$ be a finitely generated $\mathbb{N}$-graded commutative algebra with homogeneous generators $\{f_1,\cdots,f_n\}$. Then there is a positive integer $d$, such that the $d$-Veronese subring $B$ of $A$, that is, $B$ is a graded algebra with $B_i=A_{di}$, is generated in degree $1$ over $B_0$ by \cite[Exercise 7.4.G]{Vak}. Let $X=\Proj A$ and $Y=\Proj B$. Then by \cite[Exercise 7.4.D]{Vak}, $X\cong Y$. By Serre's theorem,
    \[(\coh X,\cO_X,-\otimes\cO_x(d))\cong (\coh Y,\cO_Y,-\otimes\cO_Y(1))\cong (\qgr B,\cB,s).\]

    This shows for any projective scheme $X$ with coordinate ring $A$, $\coh X$ is equivalent to a quotient category $\qgr B$ for some Veronese subring $B$ of $A$ where $B$ is generated in degree $1$ over $B_0$. The projective scheme associating to $A$ and $B$ are isomorphic. But the noncommutative projective schemes associating to $A$ and $B$ are not isomorphic in general.
\end{remark}

\subsection{Noncommutative projective coordinate rings}
Based on Corollary \ref{coor. ring}, we give a definition of coordinate rings of noncommutative projective schemes.

\begin{definition} For a noncommutative projective scheme $(\cC,O,s)$, $B(\cC,O,s)$ is called the \textit{noncommutative projective coordinate ring} of $(\cC,O,s)$.
\end{definition}

A noncommutative projective coordinate ring is exactly a right noetherian commonly graded algebra $A$ such that $A$ satisfies $\chi_1$ and $\depth_A(A)\geqslant 2$ by Corollary \ref{coor. ring}. In this case, $(\qgr A,\cA,s)$ is a noncommutative projective scheme, and $A \cong B(\qgr A,\cA,s)$.

If $A$ is a noetherian commonly graded AS-Gorenstein algebra of dimension $d\geqslant 2$, then $\depth_A(A)\geqslant 2$ by definition and $A$ satisfies $\chi$ by Theorem \ref{AS-G has balanced dualizing complex and chi}. So $A$ is a noncommutative projective coordinate ring.
In particular, AS-regular algebras are regarded as the coordinate rings of noncommutative projective spaces. So we give the following definition.

\begin{definition}\label{noncommutative-proj-space} A \textit{noncommutative projective space} is a noncommutative quasi-projective space equivalent to $(\qgr A,s)$  for some right noetherian commonly graded AS-regular algebra $A$ of dimension $d\geqslant 2$.
\end{definition}

\begin{remark}
    For a noncommutative quasi-projective space $(\cC,s)$, after fixing a structure sheaf $O$ and taking $A=B(\cC,O,s)$, $(\cC,s)$ is equivalent to the noncommutative quasi-projective space $(\qgr A,s)$, and $(\cC,O,s)$ is isomorphic to the noncommutative projective scheme associated with $A$. So, to study  noncommutative quasi-projective spaces and noncommutative projective schemes, it is enough to study $(\qgr A,s)$ and $(\qgr A,\cA,s)$ where $A$ is a right noetherian commonly graded algebra $A$ such that $A$ satisfies $\chi_1$ and $\depth_A(A)\geqslant 2$.
\end{remark}

Next we study the structure sheaves of a noncommutative quasi-projective space.

\begin{lemma}\label{oemga pi M is finitely generated}
Let $A$ be a right noetherian commonly graded algebra satisfying $\chi_1$. If $(\qgr A,\cM,s)$ is a noncommutative projective scheme, then $\omega\!\cM$ is a finitely generated graded $A$-module.
\end{lemma}
\begin{proof}
Since $(\cM,s)$ is ample in $\qgr A$ by Definition \ref{def-noncom-proj-scheme}, there are positive integers $r_1,\cdots$, $r_p$ and an epic morphism
  \[\oplus_{i=1}^ps^{-r_i}\cM\to \cA.\]
So there is an injective morphism
  \[\gHom_{\cA}(\cA,\cM)\to \gHom_{\cA}(\cM,\oplus_{i=1}^ps^{-r_i}\cM).\]
Since $\gHom_{\cA}(\cM,\cM)$ is bounded-below, it follows that $\gHom_{\cA}(\cM,\oplus_{i=1}^ps^{-r_i}\cM)$ is also bounded-below. Then $\omega\cM\cong \gHom_A(A,\omega \cM)\cong \gHom_{\cA}(\cA,\cM)$ is bounded-below. By Proposition \ref{facts about chi condition}, $(\omega\!\cM)_{\geqslant n}$ is contained in a finitely generated graded submodule of $\omega\!\cM$  for $n \in \mathbb{Z}$. Consequently, $\omega\!\cM$ is finitely generated.
\end{proof}

\begin{proposition}\label{one to one: sheaf of space and module when coor. ring fixed}
Let $A$ be a noncommutative projective coordinate ring.
Then the iso-classes of the structure sheaves $\cM$ of $(\qgr A,s)$ are in one-to-one correspondence to the iso-classes of the finitely generated graded $A$-modules $N$ such that
\begin{enumerate}
    \item [(1)] $\depth_A(N)\geqslant 2$, and
    \item [(2)] $(\cN,s)$ is ample in $\qgr A$.
\end{enumerate}
The correspondence is given by $[\cM] \mapsto [\omega\!\cM]$.
\end{proposition}
\begin{proof}
Suppose that $\cM$ is a structure sheaf of $(\qgr A,s)$. Then $\omega\!\cM$ is finitely generated by Lemma \ref{oemga pi M is finitely generated}. By Corollary \ref{depth geqslant 2}, $\depth_A(\omega\!\cM)\geqslant 2$. Since $\cM\cong \pi\omega\!\cM$, $(\pi\omega\!\cM,s)$ is ample in $\qgr A$.

Conversely, suppose $N$ is a finitely generated graded $A$-module  such that $(\cN,s)$ is ample in $\qgr A$ and $\depth_A(N)\geqslant 2$. Then $\cN$ is a noetherian object in $\qgr A$.
Since $\depth_A(N)\geqslant 2$, $N\cong \omega \cN$ by Corollary \ref{depth geqslant 2}. Hence $\gHom_{\cA}(\cN,\cN)\cong \gHom_A(N,N)$ is commonly graded.
Since $A$ satisfies $\chi_1$, $(\omega\cX)_n$ is finite-dimensional for any $X\in \gr A$ and $n\in \mathbb{Z}$ by Proposition \ref{facts about chi condition} (2). Hence for any $\cX\in \qgr A$, $\Hom_{\cA}(\cN,\cX)\cong \Hom_{\Gr A}(N,\omega\cX)$ is finite-dimensional.
Thus $(\qgr A,\cN,s)$ is a noncommutative projective scheme and $\cN$ is a structure sheaf of $(\qgr A,s)$.

At last, the one-to-one correspondence follows from $\pi\omega \cM\cong \cM$ and $N\cong \omega\pi N$ when $\depth_A(N)\geqslant 2$.
\end{proof}

\section{Equivalences between quasi-projective spaces}\label{Morita Theory for Noncommutative quasi-projective spaces}
In this section we investigate equivalences between the noncommutative quasi-projective spaces in preparation for the study of noncommutative resolutions. We prove that such equivalences are induced by modulo-torsion-invertible bimodules (see Definition \ref{definition of modulo-torsion-invertible}, and Theorems \ref{morita theory}, \ref{modulo-torsion-invertible bimodule and equivalent functor}).

\subsection{Morita-like theory of noncommutative quasi-projective spaces}
 Here are some preparation lemmas. The following lemma is a commonly graded version of \cite[Proposition 2.5(2)]{AZ}.
\begin{lemma}\label{tails A cong tails B induced by morphism}
Let $\varphi:A\to B$ be a morphism of commonly graded algebras. Suppose both the kernel and cokernel of $\varphi$ are bounded-above.
\begin{enumerate}
    \item [(1)] If the algebra $B$ is finitely generated, then so is $A$.
    \item [(2)] If  $J_A/J_A^2$ is finite-dimensional, then $-\otimes_A B:\Gr A\to \Gr B$ induces an equivalent functor  $\Qgr A \to \Qgr B$, with a quasi-inverse induced by the  restriction functor $\gHom_B({}_AB_B,-):\Gr B\to \Gr A$.
    \item [(3)] If $A$ and $B$ are right noetherian, then $-\otimes_AB$ induces an equivalent functor  $\qgr A \to \qgr B$.
\end{enumerate}
\end{lemma}
\begin{proof}
    (1) Since both $\Ker \varphi$ and $\Coker \varphi$ are bounded-above, there is an integer $n$ such that $A_{\geqslant n+1}\cong B_{\geqslant n+1}$. We may assume that  $B_{\leqslant n}$ contains a set of  generators of the algebra  $B$.
    Then, for any homogeneous element $b\in B$, $b=\sum_i b_i$ where $b_i$ is a product of homogeneous elements with degree no more than $n$.

    Let $n'=3n$ and $\deg b \geqslant n'+1$.
    Without loss of generality, we may assume that $b=y_1y_2 \cdots y_s$ is a monomial such that each $y_j \in  B_{\leqslant n}.$
    We claim that the product $b=y_1\cdots y_s$ can be divided into several parts, say,
    $y_1\cdots y_{r_1}$, $y_{r_1+1}\cdots y_{r_2}$, $\cdots$, $y_{r_t+1}\cdots y_s$, so that the degree of each part is in $[n+1,n']$. In fact, there exist $x_1=y_1\cdots y_{r_1} \in B_{\geqslant n+1} \cap B_{\leqslant 2n}$ and  $b' \in B_{\geqslant n+1}$ such that $b=x_1b'$ as $\deg b > n'=3n$. If $\deg b' \leqslant n'$, then we are done. Otherwise, we do the same for $b'$. In this way, $b$ can be written as a product of some elements with degrees in $[n+1,n']$.

    It follows from $A_{\geqslant n+1}\cong B_{\geqslant n+1}$ that for every homogeneous element $a\in A_{\geqslant n'+1}$, $a$ can be written as a finite sum of products of finite homogeneous elements in $A_{n+1}\oplus \cdots \oplus A_{n'}$.
    Therefore $A_{\leqslant n'}$ contains a set of generators of the algebra $A$,
    and $A$ is a finitely generated commonly graded algebra.

    (2) Both $\Ker \varphi$ and $\Coker \varphi$  are finite-dimensional as they are bounded-above. Thus $A/\Ker\varphi$ is a finitely presented graded left $A$-module.

    Since $J_A/J_A^2$ is finite-dimensional, $A/J_A$ is finitely presented as a graded left $A$-module by Lemma \ref{finitely generated algebra}. Thus every graded left simple $A$-module is finitely presented. Since $\Coker\varphi$ is finite-dimensional, it is finitely presented as a graded left $A$-module by an induction on the length of $\Coker\varphi$.

    It follows from the exact sequence
    \[0\to A/\Ker\varphi \to B\to \Coker\varphi\to 0\]
    that $B$ is a finitely presented graded left $A$-module. Thus there is a finitely generated minimal projective presentation $P_1\to P_0\to {}_AB\to 0$. Then for any graded simple $A$-module $X$, $\Tor_1^A(X,B)\cong X\otimes_A P_1$ which is finite-dimensional. Since ${}_AB$ is finitely generated, $X\otimes_AB$ is finite-dimensional. Hence for any finite-dimensional graded $A$-module $Y$, $Y\otimes_A B$ and $\Tor_1^A(Y,B)$ are finite-dimensional by an induction on the length of $Y$. Since $-\otimes_AB$ and $\Tor_1^A(-,B)$ commute with direct limits, $M\otimes_AB$ and $\Tor_1^A(M,B)$ are torsion for any torsion $A$-module $M$.
    By Proposition \ref{induced functors on quotient categories}, $-\otimes_A B$ induces a functor
    \[-\otimes_{\cA}\cB:\Qgr A\to \Qgr B.\]

    On the other hand, the restriction of any torsion $B$-module $K$ to $A$ is $A$-torsion, that is, $\gHom_B({}_AB_B,K)$ is a torsion $A$-module. 
    Obviously, $\gExt_B^1({}_AB_B,-)=0$. So by Proposition \ref{induced functors on quotient categories}, $\gHom_B({}_AB_B,-)$ induces a functor
    \[\gHom_{\cB}(\cB,-):\Qgr B\to \Qgr A.\]

    For any graded $A$-module $M$, by applying $M\otimes_A -$ to
    \[0\to \Ker\varphi\to A\to B\to \Coker\varphi\to 0,\]
    it follows from \cite[Proposition 2.4(5)]{AZ} that both the kernel and cokernel of $M\to M\otimes_AB$ are torsion $A$-modules. Thus $\pi_A M\to \pi_A(M\otimes_AB)=\pi_A\gHom_B({}_AB_B, M\otimes_AB)$ is an isomorphism.

    For any graded $B$-module $N$, restricting $N$ to $A$, then both the kernel and  cokernel of the morphism $N\to N\otimes_A B$ are torsion $A$-modules. In fact $N\to N\otimes_A B$ is also a $B$-module morphism. Hence both the  kernel and cokernel are also torsion $B$-modules. Thus $\pi_B N\to \pi_B(N\otimes_AB)$ is an isomorphism.
    Since the composition of $N\to N\otimes_A B\to N$ is the identity map, where $N\otimes_A B\to N$ is the canonical map, $\pi_B N\to \pi_B(N\otimes_A B)$ is the inverse of $\pi_B(N\otimes_AB) \to \pi_BN$. So $\pi_B (N\otimes_AB)=\pi_B(\gHom_B({}_AB_B,N)\otimes_AB)\to \pi_B N$ is an isomorphism.

    Consequently, $-\otimes_AB$ induces an equivalent functor on $\Qgr A\to \Qgr B$ with the quasi-inverse induced by $\gHom_B({}_AB_B,-)$.

   (3) If both $A$ and $B$ are noetherian, then $J_A/J_A^2$ is finite-dimensional by Lemma \ref{finitely generated algebra}. The proof of the noetherian case is similar.
\end{proof}

\begin{lemma}\label{lemma: equi-func induced by Hom}
Let $A$ be a right noetherian commonly graded algebra satisfying
$\chi_1$, $M$  a finitely generated graded $A$-module. Let $B=\gEnd_{\cA}(\cM)$ and $B'=\gEnd_A(M)$.

If $(\qgr A,\cM,s)$ is a noncommutative projective scheme, then \[\pi_B\gHom_{\cA}(\cM,-):\qgr A\to \qgr B\]
is induced by
\[\gHom_A(M,-)\otimes_{B'}B:\gr A\to \Gr B.\]
\end{lemma}
\begin{proof}
By Theorem \ref{noncommutative Serre theorem for commonly graded algebra}, $B$ is a right noetherian commonly graded algebra and there is an isomorphism of noncommutative projective schemes
\[\pi_B\gHom_{\cA}(\cM,-):(\qgr A,\cM,s)\to (\qgr B,\cB,s), \,\, \cN\mapsto \pi_B\gHom_{\cA}(\cM,\cN).\]

Let $0\to K\to P\to M\to 0$ be an exact sequence of finitely generated graded $A$-modules where $P$ is projective. For any $N\in \tor A$, $\gHom_A(P,N)$ is finite-dimensional. It follows from the exact sequence
\begin{equation}\label{equ.of-ext-hom(M,N)}
0\to \gHom_A(M,N)\to \gHom_A(P,N) \to \gHom_A(K,N)\to \gExt_A^1(M,N)\to 0 \tag{$\star$}
\end{equation}
that $\gHom_A(M,N)$ is finite-dimensional. So $\gHom_A(M,N)\in \tor B'$.

Similarly, we have $\gHom_A(K,N)$ is finite-dimensional. Hence by the exact sequence \eqref{equ.of-ext-hom(M,N)}, $\gExt_A^1(M,N)$ is finite-dimensional. By Proposition \ref{induced functors on quotient categories}, $\gHom_A(M,-)$ induces a functor $F_1:\qgr A\to \Qgr B'$.

Since $A$ satisfies $\chi_1$, the natural map $B'=\gHom_A(M,M)\to B=\gHom_{\cA}(\cM,\cM)$ has bounded-above kernel and cokernel by Proposition \ref{facts about chi condition} (3). Since $B$ is a right noetherian commonly graded algebra, it is a finitely generated algebra by Lemma \ref{finitely generated algebra}. By Lemma \ref{tails A cong tails B induced by morphism}, $-\otimes_{B'}B$ induces an equivalent functor $F_2:\Qgr B'\to \Qgr B$. So we have the following commutative diagram:
\[
\xymatrix{
\gr A \ar[r]^{\gHom_A(M,-)} \ar[d]^{\pi_A} & \Gr B' \ar[r]^{-\otimes_{B'}B} \ar[d]^{\pi_{B'}} & \Gr B \ar[d]^{\pi_B}\\
\qgr A \ar[r]^{F_1} & \Qgr B' \ar[r]^{F_2} & \Qgr B.
}\]

Let $F=F_2F_1$. Then $F$ is induced by the functor $\gHom_A(M,-)\otimes_{B'}B:\gr A\to \Gr B$. We claim that $F\cong \pi\gHom_{\cA}(\cM,-)$. Then $F$ can be regarded as a functor from $\qgr A$ to $\qgr B$ by Theorem \ref{noncommutative Serre theorem for commonly graded algebra} and the proof is completed.

By Proposition \ref{facts about chi condition}, for any $N\in \gr A$, the natural map $\gHom_A(M,N)\to \gHom_{\cA}(\cM,\cN)$ has bounded-above kernel and cokernel. Regarding $\gHom_{\cA}(\cM,\cN)$ as a graded $B'$-module via the restriction functor, then
\[\pi_{B'}\gHom_A(M,N)\cong \pi_{B'} \gHom_{\cA}(\cM,\cN)\]
in $\Qgr B'$.
By Lemma \ref{tails A cong tails B induced by morphism}, the inverse of $F_2$ is induced by the restriction functor. Thus $F_2(\pi_{B'} \gHom_{\cA}(\cM,\cN))\cong \pi_B\gHom_{\cA}(\cM,\cN)$.
Hence there is a natural isomorphism
\[F(\cN)\cong F_2(\pi_{B'}\gHom_A(M,N))\cong F_2(\pi_{B'} \gHom_{\cA}(\cM,\cN))\cong \pi_B\gHom_{\cA}(\cM,\cN).
\qedhere\]
\end{proof}

If $\depth_A(M)\geqslant 2$, then $\gEnd_A(M)\cong \gEnd_{\cA}(\cM)$ as graded algebras by Corollary \ref{depth geqslant 2}. Thus we have the following corollary.

\begin{corollary}\label{corollary: equi-func induced by Hom}
Let $A$ be a right noetherian commonly graded algebra satisfying $\chi_1$, and $M$ be a finitely generated graded $A$-module with $\depth_A(M)\geqslant 2$. Let $B=\gEnd_{\cA}(\cM)$.
If $(\qgr A,\cM,s)$ is a noncommutative projective scheme, then $\pi_B\gHom_{\cA}(\cM,-):\qgr A\to \qgr B$ is induced by $\gHom_A(M,-):\gr A\to \Gr B$.
\end{corollary}

We are now ready for studying Morita equivalences between the noncommutative quasi-projective spaces.

\begin{proposition}\label{M to (B,F)}
Let $A$ be a right noetherian commonly graded algebra satisfying $\chi_1$.
Consider the following two statements.
    \begin{enumerate}
        \item [(1)] $M$ is a finitely generated graded $A$-module with $\depth_A(M)\geqslant 2$, and $(\qgr A,\cM,s)$ is a noncommutative projective scheme.
        \item [(2)] $F: \qgr A\to \qgr B$ is an equivalent functor commuting with the polarizations, where $B$ is a noncommutative projective coordinate ring.
    \end{enumerate}
Then, the graded $A$-modules satisfying (1) are in one-to-one correspondence to the pairs $(B,F)$ satisfying (2) up to isomorphism. The correspondence is given by
\[M\mapsto (\gEnd_A(M),\pi_{\gEnd_A(M)}\gHom_{\cA}(\cM,-)),\]
with the inverse given by
$(B,F)\to \omega F^{-1}\cB.$
\end{proposition}
\begin{proof}
Suppose that $M$ is a graded $A$-module satisfying (1).
It follows from $\depth_A(M)\geqslant 2$ that $\gEnd_A(M)\cong \gEnd_{\cA}(\cM)$ by Corollary \ref{depth geqslant 2}.
Let $B=\gEnd_A(M)$. Then by Proposition \ref{coor. ring}, $B$ is a right noetherian commonly graded algebra satisfying $\chi_1$ and $\depth_B(B)\geqslant 2$. By Theorem \ref{noncommutative Serre theorem for commonly graded algebra}, $F=\pi_B\gHom_{\cA}(\cM,-):\qgr A\to \qgr B$  is an equivalent functor which gives an isomorphism between the noncommutative projective schemes $(\qgr A,\cM,s)\to (\qgr B,\cB,s)$.

Suppose that $(B,F)$ is a pair satisfying (2).
Let $G$ be a quasi-inverse of $F$ and $M= \omega G\cB$. Then $\depth_A(M)\geqslant 2$ by Corollary \ref{depth geqslant 2} and $\cM\cong G\cB$. Since $F:\qgr A\to \qgr B$ is an equivalent functor commuting with $s$, $(\qgr A,\cM,s)$ is a noncommutative projective scheme. By Lemma \ref{oemga pi M is finitely generated}, $M_A$ is finitely generated.

To prove the correspondence is one-to-one, it remains to show that $B\cong \gEnd_A(M)$ and $F\cong \pi_{\gEnd_A(M)}\gHom_{\cA}(\cM,-))$.

Since $\depth_B(B)\geqslant 2$ and $\depth_A(M)\geqslant 2$, it follows from Corollary \ref{depth geqslant 2} that
\[B\cong \gHom_B(B,B)\cong \gHom_{\cB}(\cB,\cB)\cong \gHom_{\cA}(\cM,\cM)\cong \gHom_A(M,M).\]
For any $\cX\in \qgr A$,
\begin{align*}
    F\cX &\cong \pi_B\omega_B F\cX\\
         &\cong \pi_B\gHom_B(B,\omega_BF\cX)\\
         &\cong \pi_B\gHom_{\cB}(\cB,F\cX)\\
         &\cong \pi_B\gHom_{\cA}(\cM,\cX).
\end{align*}
Hence $F\cong \pi_B\gHom_{\cA}(\cM,-):\qgr A\to \qgr B$.
\end{proof}

The following theorem is the main result of this section.

\begin{theorem}\label{morita theory}
Let $A$ and $B$ be two noncommutative projective coordinate rings.
Suppose
\[F:(\qgr A,s)\rightleftarrows (\qgr B,s): G\]
is an equivalence of the noncommutative quasi-projective spaces. Let $M=\omega G\cB$ and $M'=\omega F\cA$. Then
\begin{enumerate}
    \item [(1)] $M_A$ and $M'_B$ are finitely generated.
    \item [(2)] $F\cong \pi_B\gHom_{\cA}(\cM,-)$ and $G\cong\pi_A\gHom_{\cB}(\cM',-)$.
    \item [(3)] $B\cong \gEnd_A(M)$ and $A\cong \gEnd_B(M')$.
    \item [(4)] $M\cong\gHom_B(M',B)$ and $M'\cong \gHom_A(M,A)$ as graded bimodules.
\end{enumerate}
     Hence there is a graded Morita context $(A,B,M,M',\tau,\mu)$ isomorphic to the graded Morita context defined by $M_A$ or $M'_B$.

\begin{enumerate}
    \item [(5)] Suppose $A$ (resp. $B$) is also left noetherian. Then ${}_AM'$ (resp. ${}_BM$) is finitely generated, and $F\cong -\otimes_{\cA}\cM'$ (resp. $G\cong -\otimes_{\cB}\cM$).
\end{enumerate}

\begin{enumerate}
    \item [(6)] Suppose both $A$ and $B$ are left noetherian. Then $F,G$ are induced by the graded Morita context $(A,B,M,M',\tau,\mu)$.
    \item [(7)] Suppose both $A$ and $B$ are left noetherian. Then $(A,B,M,M',\tau,\mu)$ induces an equivalence
    \[\cM\otimes_{\cA}-:\qgr A^o\rightleftarrows \qgr B^o: \cM'\otimes_{\cB}-.\]
\end{enumerate}
\end{theorem}
\begin{proof}
Note that $(\cM,s)$ is ample in $\qgr A$ and $(\cM',s)$ is ample in $\qgr B$.

(1), (2) and (3) follow from Proposition \ref{M to (B,F)}.

(4) Since $\depth_A(M)\geqslant 2$, $M\cong \omega\cM$ by Corollary \ref{depth geqslant 2}. Then $M\cong \omega G\cB$. Since $G\cong \pi_A\gHom_{\cB}(\cM',-)$, it follows that $M\cong \omega\pi_A\gHom_{\cB}(\cM',\cB)$. 
By Corollary \ref{depth of coordinate ring of noncom. proj. scheme}, there is an isomorphism $\omega\pi_A\gHom_{\cB}(\cM',\cB)\cong \gHom_{\cB}(\cM',\cB)$. Since $\depth_B(B)\geqslant 2$, it follows that  $B\cong \omega\cB$ and $M\cong \gHom_{\cB}(\cM',\cB)\cong \gHom_B(M',B)$.

Similarly $M'\cong \gHom_A(M,A)$. It follows from  the definition that there is a graded Morita context $(A,B,M,M',\tau,\mu)$ isomorphic to the graded Morita context defined by $M_A$ or $M'_B$.

(5) Let $P_A\to M_A$ be a surjective morphism of graded $A$-modules where $P_A$ is a finitely generated free $A$-module. Then it induces an injective morphism of graded $A^o$-modules $M'\to \gHom_A(P,A)$. Since $A$ is left noetherian, 
${}_AM'$ is finitely generated.
We may assume $B=\gEnd_A(M)$ and $M'=\gHom_A(M,A)$ by (3) and (4). By Lemma \ref{induced functor on tails by tensor}, $M'$ is compatible with $\tor A$ and $\tor B$. Thus $-\otimes_AM'$ induces a functor
\[-\otimes_{\cA}\cM':\qgr A\to \qgr B.\]

Next we prove $F\cong -\otimes_{\cA}\cM'$. It suffices to prove $\pi_B\gHom_{\cA}(\cM,-)\cong -\otimes_{\cA}\cM'$.
By Corollary \ref{corollary: equi-func induced by Hom}, $\pi_B\gHom_{\cA}(\cM,-)$ is induced by $\gHom_A(M,-)$. So, we only need to show $\pi_B\gHom_A(M,-)\cong \pi_B(-\otimes_AM')$ by Proposition \ref{induced functors on quotient categories}.

For any $N\in \gr A$, there is a natural morphism
\[\varphi_N:N\otimes_AM'\to\gHom_A(M,N)\]
such that $\varphi_N(n\otimes f)(m)=nf(m)$ for any $n\otimes f\in N\otimes_AM'$ and $m \in M$.

Let $0\to K\to P\to N\to 0$ be an exact sequence of finitely generated graded $A$-modules with $P_A$ projective. Then there is a commutative diagram
\[
\xymatrix{
 & K\otimes_A M' \ar[r] \ar[d]^{\varphi_K} & P\otimes_A M' \ar[r] \ar[d]^{\varphi_P} & N\otimes_A M' \ar[r] \ar[d]^{\varphi_N} & 0\\
 0 \ar[r] & \gHom_A(M,K) \ar[r] & \gHom_A(M,P) \ar[r] & \gHom_A(M,N). &
}\]
Since $P_A$ is a finitely generated projective module, $\varphi_P$ is an isomorphism.

Since $\pi_AP\to \pi_AN$ is epic and $(\cM,s)$ is ample in $\qgr A$, there is an integer $n_0$ such that the natural map
\[\gHom_{\cA}(\cM,\pi_A P)_{\geqslant n_0}\to \gHom_{\cA}(\cM,\pi_A N)_{\geqslant n_0}\]
is surjective. It follows from Proposition \ref{facts about chi condition} that there are integers $n_1,n_2$ such that the natural maps
\[\gHom_A(M,P)_{\geqslant n_1}\to \gHom_{\cA}(\cM,\pi_AP)_{\geqslant n_1},\]
\[\gHom_A(M,N)_{\geqslant n_2}\to \gHom_{\cA}(\cM,\pi_AN)_{\geqslant n_2}\]
are bijective.
Let $n=\max\{n_0,n_1,n_2\}$. It follows that the natural map
\[\gHom_A(M,P)_{\geqslant n}\to \gHom_A(M,N)_{\geqslant n}\]
is surjective. Hence $\Coker\varphi_N$ is bounded-above.
Obviously, $\Ker \varphi_N\cong \Coker\varphi_K$. As $N$ is arbitrary, $\Coker\varphi_K$ is bounded-above. It follows that $\Ker\varphi_N$ is bounded-above. Consequently,
\[\pi_B(\varphi_N):\pi_B(N\otimes_AM')\to \pi_B\gHom_A(M,N)\]
is an isomorphism. So $\pi_B\gHom_A(M,-)\cong \pi_B(-\otimes_AM')$ and $F\cong -\otimes_{\cA}\cM'$.

(6) It follows from (5).

(7) By Theorem \ref{equ. of quot. cat. induced by Morita context for commonly graded}, both $\Coker\tau$ and $\Coker \mu$ are finite-dimensional. Then the assertion follows from a left module version of Theorem \ref{equ. of quot. cat. induced by Morita context for commonly graded}.
\end{proof}

\begin{corollary}\label{center of noncommutative coordinate ring}
    If $A$ and $B$ are noetherian noncommutative projective coordinate rings such that $(\qgr A,s)\cong (\qgr B,s)$, then the center of $A$ is isomorphic the center of $B$.
\end{corollary}
\begin{proof}
    By Theorem \ref{morita theory}, there is a graded $(B,A)$-bimodule $M$ such that $\gEnd_A(M)\cong B$ and $\gEnd_{B^o}(M)\cong A^o$. Hence the center of $A$ is isomorphic the center of $B$.
\end{proof}

\subsection{Equivalences and modulo-torsion-invertible bimodules}
Now we study the modulo-torsion-invertible bimodules (Definition \ref{definition of modulo-torsion-invertible}) of noncommutative projective coordinate rings.

The following lemma is direct from Theorem \ref{equ. of quot. cat. induced by Morita context for commonly graded}.
\begin{lemma}
    Let $A$ and $B$ be noetherian noncommutative projective coordinate rings. If ${}_BM_A$ is a modulo-torsion-invertible bimodule with the inverse $M'$ associated to the graded Morita context $(A,B,M,M',\tau,\mu)$.
    Then $(A,B,M,M',\tau,\mu)$ induces an equivalence of noncommutative quasi-projective spaces
         \[-\otimes_{\cA}\cM':(\qgr A,s)\rightleftarrows (\qgr B,s): -\otimes_{\cB}\cM.\]
\end{lemma}

\begin{proposition}\label{modulo-torsion-invertible of coordinate ring}
    Let $A$ and $B$ be noetherian noncommutative projective coordinate rings. If ${}_BM_A$ is a modulo-torsion-invertible bimodule with the inverse $M'$ associated to the graded Morita context $(A,B,M,M',\tau,\mu)$. Let $N=\omega_A\pi_AM$ and $N'=\omega_B\pi_BM'$. Then
    \begin{enumerate}
        \item [(1)] ${}_BN_A$ and ${}_AN'_B$ are finitely generated on both sides.
        \item [(2)] $\depth_AN\geqslant 2$ and $\depth_BN'\geqslant 2$.
        \item [(3)] ${}_BN_A$ is a modulo-torsion-invertible bimodule with the inverse ${}_AN'_B$ associated to the graded Morita context defined by ${}_BN_A$.
        \item [(4)] The equivalent functors $-\otimes_{\cA}\cM'$ and $-\otimes_{\cA}\cN':\qgr A\to \qgr B$, induced by $-\otimes_AM'$ and $-\otimes_AN'$ respectively, are naturally isomorphic.
        \item [(5)] The equivalent functors $-\otimes_{\cB}\cM$ and $-\otimes_{\cB}\cN:\qgr B\to \qgr A$, induced by $-\otimes_BM$ and $-\otimes_BN$ respectively, are naturally isomorphic.
    \end{enumerate}
    In particular, if $\depth_AM\geqslant 2$ (resp. $\depth_BM'\geqslant 2$), then the graded Morita context defined by $M_A$ (resp. $M'_B$) and the graded Morita context $(A,B,M,M',\tau,\mu)$ induce the same equivalence between $(\qgr A,s)$ and $(\qgr B,s)$ up to isomorphism.
\end{proposition}
\begin{proof}
   (1), (2) and (3) follow from Corollary \ref{depth geqslant 2}, Theorem \ref{equ. of quot. cat. induced by Morita context for commonly graded} and Theorem \ref{morita theory} because $N=\omega_A\pi_A M\cong \omega_A(\cB\otimes_{\cB}\cM)$ and $N'=\omega_B\pi_B M'\cong \omega_B(\cA\otimes_{\cA}\cM')$.

    (4) Note that the inverses of $-\otimes_{\cA}\cM'$ and $-\otimes_{\cA}\cN'$ are $-\otimes_{\cB}\cM$ and $-\otimes_{\cB}\cN$ respectively. Since $N\cong \omega_A(\cB\otimes_{\cB}\cM) \cong \omega_A(\cB\otimes_{\cB}\cN)$, (4) follows from Proposition \ref{M to (B,F)}.

    (5) Similar to (4).
\end{proof}

\begin{remark}
    By Proposition \ref{modulo-torsion-invertible of coordinate ring}, for noetherian noncommutative coordinate rings $A$ and $B$, every modulo-torsion-invertible bimodule ${}_BM_A$, associated to the graded Morita context $(A,B,M,M',\tau,\mu)$, can be replaced by a module-torsion-invertible bimodule ${}_BN_A$ with $\depth_AN\geqslant 2$, the inverse ${}_AM'_B$ of ${}_BM_A$ can be replaced by $N'=\gHom_A(N,A)$, and the  Morita context $(A,B,M,M',\tau,\mu)$ can be replaced by the graded Morita context defined by $N_A$, as the graded Morita context defined by $N_A$ and the graded Morita context $(A,B,M,M',\tau,\mu)$ induce the same equivalence.

    So in the rest we will only consider modulo-torsion-invertible $(B,A)$-bimodules ${}_BM_A$  associated to the graded Morita context defined by $M_A$ with $\depth_AM\geqslant 2$. If there is no ambiguity, we will always omit this graded Morita context defined by $M_A$ when $\depth_AM\geqslant 2$.
\end{remark}

There are some characterizations of modulo-torsion-invertible bimodules following Theorem \ref{morita theory}.
\begin{theorem}\label{modulo-torsion-invertible bimodule and equivalent functor}
Let $A$ and $B$ be noetherian noncommutative projective coordinate rings. 
\begin{enumerate}
    \item [(1)] If $M$ is a graded $(B,A)$-bimodule with $\depth_AM\geqslant 2$, then $M$ is a modulo-torsion-invertible $(B,A)$-bimodule if and only if $M$ is a finitely generated graded $A$-module such that $(\cM,s)$ is ample in $\qgr A$ and $\gEnd_A(M)\cong B$.
    \item [(2)] Every equivalence between noncommutative quasi-projective spaces $(\qgr A,s)$ and $(\qgr B,s)$ is induced by the Morita context defined of a modulo-torsion-invertible $(B,A)$-bimodule $M$ with $\depth_AM\geqslant 2$.
    \item [(3)] 
    There is a one-to-one correspondence between the equivalences $F:(\qgr A,s)\to (\qgr B,s)$ of the noncommutative quasi-projective spaces and the modulo-torsion-invertible  $(B,A)$-bimodules $M$ with $\depth_AM\geqslant 2$ up to isomorphisms.
    The correspondence is given by
\[F\mapsto \omega F^{-1}\cB.\]
The inverse is given by
\[M\mapsto -\otimes_{\cA}\cM'\cong\pi_B\gHom_{\cA}(\cM,-)\]
where $M'=\gHom_A(M,A)$.
\end{enumerate}
\end{theorem}
\begin{proof}
It follows from Proposition \ref{M to (B,F)} and Theorem \ref{morita theory}.
\end{proof}

Next proposition gives the relation between the composition of equivalent functors of noncommutative quasi-projective spaces and the corresponding modulo-torsion-invertible bimodules.

\begin{proposition}
Let $A$, $B$ and $C$ be noetherian noncommutative projective coordinate rings, and
\[(\qgr A,s)\xrightarrow[]{F_1} (\qgr B,s)\xrightarrow[]{F_2}(\qgr C,s)\]
be equivalences of noncommutative quasi-projective spaces.
Suppose
\begin{enumerate}
    \item [(1)] $F_1$ corresponds to the modulo-torsion-invertible $(B,A)$-bimodule $M_1$;
    \item [(2)] $F_2$ corresponds to the modulo-torsion-invertible $(C,B)$-bimodule $M_2$.
\end{enumerate}
Let $M_1'=\gHom_A(M_1,A)$.
Then $F_2F_1$ corresponds to the modulo-torsion-invertible $(C,A)$-bimodule $\gHom_B(M_1',M_2)$.
\end{proposition}
\begin{proof}
    By Theorem \ref{morita theory},
    $F_1^{-1}\cong \pi_A\gHom_{\cB}(\cM_1',-)$ and $F_2^{-1}\cong \pi_B\gHom_{\cC}(\cM_2',-)$ where $M_1'=\gHom_A(M_1,A)$ and $M_2'= \gHom_B(M_2,B)$.
    So
    \[F_2^{-1}(\cC)\cong \pi_B\gHom_{\cC}(\cM_2',\cC)\cong \pi_B\gHom_C(M_2',C)\cong \pi_BM_2=\cM_2.\]
    Then
    \[F_1^{-1}F_2^{-1}(\cC)\cong F_1^{-1}(\cM_2)\cong \pi_A\gHom_{\cB}(\cM_1',\cM_2)\cong \pi_A\gHom_B(M_1',M_2).\]

    Since $\depth_B(M_2)\geqslant 2$, $\gHom_{\cB}(\cM_1',\cM_2)\cong \gHom_B(M_1',M_2)$ by Corollary \ref{depth geqslant 2}. It follows from Corollary \ref{depth of coordinate ring of noncom. proj. scheme} that $\depth_A(\gHom_B(M_1',M_2))\geqslant 2$. Hence
    \[\omega F_1^{-1}F_2^{-1}(\cC)\cong \gHom_B(M_1',M_2).\]
So, $F_2F_1$ corresponds to the modulo-torsion-invertible $(C,A)$-bimodule $\gHom_B(M_1',M_2)$.
\end{proof}

We end up this section by a result which will be used later.

\begin{theorem}\label{equivalence: modules with depth geqslant 2}
    Let $A,B$ be noetherian noncommutative projective coordinate rings, and $M$ a modulo-torsion-invertible $(B,A)$-bimodule with $\depth_A M\geqslant 2$.
    Then
        \[\gHom_A(M,-):\{N\in \gr A\mid \depth_A(N)\geqslant 2\}\to \{Y\in\gr B\mid \depth_B(Y)\geqslant 2\}\]
    gives an equivalence between the full subcategories of $\gr A$ and $\gr B$, with a quasi-inverse $\gHom_B(M',-)$ where $M'=\gHom_A(M,A)$.
\end{theorem}
\begin{proof}
    Suppose $N$ is a finitely generated graded $A$-module and $\depth_A(N)\geqslant 2$. Then by Corollary \ref{depth geqslant 2},
    \[\gHom_A(M,N)\cong \gHom_A(M,\omega\pi N)\cong \gHom_{\cA}(\cM,\cN).\]
    It follows that $\gHom_{\cA}(\cM,\cN)$ is bounded-below. By Lemma \ref{facts extracted from noncom. Serre thm}, $\gHom_{\cA}(\cM,\cN)$ is a finitely generated graded $B$-module and so is $\gHom_A(M,N)$. 
    It follows from Corollary \ref{depth of coordinate ring of noncom. proj. scheme} that $\depth_B(\gHom_A(M,N))\geqslant 2$.
    By Corollary \ref{corollary: equi-func induced by Hom},
    \[\pi_A\gHom_B(M',\gHom_A(M,N))\cong \pi_A\gHom_{\cB}(\cM',\pi_B\gHom_A(M,N))\cong\pi_AN.\]
    Since $\depth_B(\gHom_A(M,N))\geqslant 2$, by Corollary \ref{depth geqslant 2},
    \[\gHom_B(M',\gHom_A(M,N))\cong \gHom_{\cB}(\cM',\pi_B\gHom_A(M,N)).\]
    So the depth of $\gHom_B(M',\gHom_A(M,N))$ over $A$ is no less than $2$ by Corollary \ref{depth of coordinate ring of noncom. proj. scheme}. Consequently,
    \[\gHom_B(M',\gHom_A(M,N))\cong N.\]

    By a similar argument, for any $N'\in \gr B$ with $\depth_B(N')\geqslant 2$,
    \[\gHom_A(M,\gHom_B(M',N'))\cong N'.\]
    So $\gHom_A(M,-)$ and $\gHom_B(M',-)$ give an equivalence between
    \[\{N\in \gr A\mid \depth_A(N)\geqslant 2\}\text{ and } \{Y\in\gr B\mid \depth_B(Y)\geqslant 2\}.
    \qedhere\]
\end{proof}

\section{Noncommutative isolated singularities}\label{NC-singu}
In this section, we define and study noncommutative resolutions for commonly graded AS-Gorenstein isolated singularities.

\subsection{Noncommutative isolated singularity}
Noncommutative isolated singularities, as an analogue of isolated singularities in commutative case are studied in \cite{Jo, Ue1, Ue2, MU} for $\mathbb{N}$-graded algebras.

\begin{definition}\label{def-nc-iso-sing}
  Let $A$ be a right noetherian commonly graded algebra and $D(\qgr A)$ be the derived category of $\qgr A$. For any $\cM,\cN\in \qgr A$, let
  \[\Ext_{\qgr A}^i(\cM,\cN)=\Hom_{\D(\qgr A)}(\cM,\cN[i]).\]

  The \textit{global dimension} $\gldim(\qgr A)$ of $\qgr A$ is defined to be
  \[\sup\{i\mid \Ext_{\qgr A}^i(\cM,\cN)\neq 0, \cM,\cN\in \qgr A\}.\]
  If $\gldim (\qgr A)$ is finite, then $A$ is called a \textit{noncommutative isolated singularity}.
\end{definition}

Next proposition justifies the name noncommutative isolated singularity. Recall that a commutative graded ring $R$ is called a \textit{graded isolated singularity} if  the degree $0$ part of the homogeneous localization $R_{(\mathfrak{p})}$ is a regular local ring for any non-maximal graded prime ideal $\mathfrak{p}$, that is, the associated projective scheme $\Proj R$ is smooth \cite{Ha}.

\begin{proposition}\cite[Corollary 4.13]{LW2}
Let $R=k[x_1,\cdots,x_n]/I$ be a graded quotient of the polynomial algebra $k[x_1,\cdots,x_n]$ with $\deg x_i=1$ and $\mathfrak{m}=R_{>0}$. Then the following are equivalent.
\begin{enumerate}
    \item [(1)] $R$ is a graded isolated singularity.
    \item [(2)] The homogeneous localization $R_{(\mathfrak{p})}$ is a regular graded local ring for any non-maximal graded prime ideal $\mathfrak{p}$ of $R$.
    \item [(3)] The global dimension of $\qgr R$ is finite.
    \item [(4)] The global dimension of $\coh (\Proj R)$ is finite.
\end{enumerate}
\end{proposition}

For any $\cM\in \Qgr A$, the \textit{injective dimension} $\idim_{\cA}\cM$ of $\cM$ is the length of the minimal injective resolution of $\cM$ in $\Qgr A$.
Next lemma shows that the Ext-groups $\Ext_{\qgr A}^i(\cM,\cN)$, and so the global dimension of the abelian category $\qgr A$, can be computed in $\Qgr A$, which is a Grothendieck category with enough injective objects.

\begin{lemma}\label{Ext-groups-in-tails}
Let $A$ be a right noetherian commonly graded algebra.
\begin{enumerate}
\item [(1)] For any $\cM, \cN \in \qgr A$,
$\Ext^n_{\qgr A}(\cM, \cN) \cong \Ext^n_{\Qgr A}(\cM, \cN).$
\item [(2)] $\gldim(\qgr A)=\max\{\idim_{\cA}\cN\mid \cN\in \qgr A\}$.
\end{enumerate}
\end{lemma}
\begin{proof}
See \cite[Lemmas 4.1 and 4.2]{LW2}.
\end{proof}

\begin{lemma}\label{pi preserves injective hulls}
Let $A$ be a right noetherian commonly graded algebra. Then the quotient functor $\pi:\Gr A\to \Qgr A$ preserves injective hulls.
\end{lemma}
\begin{proof}
Let $I$ be a graded injective $A$-module. By \cite[Lemma 5.1]{LW},
there is a decomposition $I=I_1\oplus I_2$ where $I_1$ is torsion and $I_2$ is torsion-free. Then the lemma follows from \cite[Lemma 4.3]{LW2}.
\end{proof}

Noncommutative isolated singularities have close relations with the minimal graded injective resolutions of modules as we will see more in Theorem \ref{minimal injective resolution over isolated singularity}.

\begin{lemma}\label{Ext in tails cong Ext in gr when depth>n}
    Let $A$ be a right noetherian commonly graded algebra, $N$ be a finitely generated graded $A$-module such that $\depth_A(N)\geqslant n$.
    \begin{enumerate}
        \item [(1)] If $0\to N\to I^0\to I^1\to \cdots$ is the minimal graded injective resolution of $N_A$, then $I^i$ is torsion-free for any $i<n$.
        \item [(2)] For any $M\in\gr A$ and $0\leqslant i<n-1$, $\gExt_A^i(M,N)\cong \gExt_{\cA}^i(\cM,\cN)$.
    \end{enumerate}
\end{lemma}
\begin{proof}
    (1) Since $\depth_A(N)\geqslant n$,
    $\gExt^i(S,N)\cong\gHom_A(S,I^i)=0$ for all $i<n$, where $S=A/J_A$.  Hence $I^i$ is torsion-free for $i < n$.

    (2) By Lemma \ref{pi preserves injective hulls}, $0\to \cN\to \pi I^0\to \pi I^1\to \cdots$ is an injective resolution of $\cN$ in $\Qgr A$. So, for any $M\in\gr A$,
    \[\gExt_{\cA}^i(\cM,\cN)=\textrm{H}^i(\gHom_{\cA}(\cM,\pi I^\bullet))\cong \textrm{H}^i(\gHom_A(M,\omega\pi I^\bullet)).\]
    Since $I^i$ is torsion-free for any $i<n$, $\omega\pi I^i\cong I^i$. Thus, for any $0\leqslant i<n-1$,
    \[\gExt_{\cA}^i(\cM,\cN)\cong \textrm{H}^i(\gHom_A(M,I^\bullet))\cong\gExt_A^i(M,N).
    \qedhere\]
\end{proof}

\begin{theorem}\label{minimal injective resolution over isolated singularity}
    Let $A$ be a right noetherian commonly graded algebra.
    \begin{enumerate}
        \item [(1)]  If $A$ is a noncommutative isolated singularity with $\gldim(\qgr A)=d-1$, then for any  finitely generated graded $A$-module $N$, $I^i$ is torsion for any $i\geqslant d$, where $I^i$ is the $i$-th term  in the minimal graded injective resolution of $N$.
        \item [(2)] Suppose that $A$ is a balanced Cohen-Macaulay algebra of dimension $d$. If for any MCM $A$-module $M$ and any $i\geqslant d$, $I^i$ is torsion where $I^i$ is the $i$-th term  in the minimal graded injective resolution of $M$,
            then $A$ is a noncommutative isolated singularity with $\gldim(\qgr A)=d-1$.
    \end{enumerate}
\end{theorem}
\begin{proof}
    (1) By Lemma \ref{pi preserves injective hulls}, $0\to \cN\to\pi I^0\to \pi I^1 \to \cdots \to \pi I^n \to \cdots $ is the minimal injective resolution of $\cN$ in $\Qgr A$.
    Since $\gldim(\qgr A)=d-1$, by Lemma \ref{Ext-groups-in-tails}, $\pi I^i=0$ for all $i\geqslant d$. Thus $I^i$ is torsion for all $i\geqslant d$.

    (2) For any finitely generated graded $A$-module $N$, let $\cdots\to P_1\xrightarrow[]{\partial_1} P_0\xrightarrow[]{\partial_0} N\to 0$ be a finitely generated graded projective resolution of $N_A$. Consider the exact sequence
    \[0\to \Ker\partial_{d-1}\to P_{d-1}\to\cdots\to P_0\to N\to 0.\]
   We claim that $\Ker\partial_{d-2}$ has a graded injective resolution in which  the $i$-th term is torsion for all $i\geqslant d$. We may assume that $\Ker\partial_{d-1}\neq 0$.

    It follows from Lemmas \ref{depth in exact sequences} and \ref{depth < infinity} that $\depth_A(\Ker\partial_{d-1})=d$.
    Hence $\Ker\partial_{d-1}$ is an MCM $A$-module.
    Let $0\to \Ker\partial_{d-1}\to E^0\to E^1\to \cdots$ and $0\to P_{d-1}\to E'^0\to E'^1\to \cdots$ be the minimal graded injective resolutions of $\Ker\partial_{d-1}$ and $P_{d-1}$ respectively. Then $E^i$ and $E'^i$ are torsion for any $i\geqslant d$. There is a chain map $\alpha:E^\bullet\to E'^\bullet$ induced by $\Ker\partial_{d-1}\to P_{d-1}$ such that the mapping cone $\cone(\alpha)$ of $\alpha$ is a graded injective resolution of $\Ker\partial_{d-2}$. Note that $\cone(\alpha)^i=E^{i+1}\oplus E'^i$. So for any $i\geqslant d$, $\cone(\alpha)^i$ is torsion.

    By induction, there is a graded injective resolution $E_N^\bullet$ of $N$ such that for any $i\geqslant d$, $E^i_N$ is torsion. It follows that the injective dimension of $\cN$ in $\Qgr A$ is at most $d-1$.

    Since $\depth_A(A)= d$, it follows from  Lemma \ref{Ext in tails cong Ext in gr when depth>n} that  the injective dimension of $\cA$ in $\Qgr A$ is $d-1$. Thus $\gldim(\qgr A)=d-1$ and so $A$ is a noncommutative isolated singularity.
\end{proof}

\begin{definition}\label{definition of GAS isolated singulrity}
    A commonly graded algebra $A$ is called a \textit{commonly graded AS-Gorenstein isolated singularity} of dimension $d$ if it is a noetherian commonly graded AS-Gorenstein algebra of dimension $d$ and it is a noncommutative isolated singularity.
\end{definition}

\begin{lemma}
If $A$ is a noetherian commonly graded AS-Gorenstein algebra of dimension $d$, then $A$ is a noncommutative isolated singularity if and only if so is  $A^o$.
\end{lemma}
\begin{proof}
    It suffices to prove one direction. Suppose $A$ is a noncommutative isolated singularity.
    Let $\underline{\MCM}(A)$ and $\underline{\MCM}(A^o)$ be the stable categories of $\MCM(A)$ and $\MCM(A^o)$ respectively, which are triangulated categories. Let $[1]$ denote the shift functors of these two triangulated categories.
    Since $\gHom_A(-,A)$ and $\gHom_{A^o}(-,A)$ send finitely generated graded projective modules to finitely generated graded projective modules, they induce a duality
    \[\gHom_A(-,A):\underline{\MCM}(A)\rightleftarrows \underline{\MCM}(A^o):\gHom_{A^o}(-,A)\]
    by Proposition \ref{Hom(-,A) induce duality on MCM}.

    Let $\Omega=D(R^d\Gamma_A(A))$. It follows from \cite[Theorem 1.3]{Ue1} that $-\otimes_A\Omega[d-1]$ is a Serre functor of $\underline{\MCM}(A)$. Then, by duality,
    \[\gHom_A(-,A)\circ(-\otimes_A\Omega[d-1])\circ \gHom_{A^o}(-,A)\]
    is a Serre functor of $\underline{\MCM}(A^o)$. Note that
    \begin{align*}
        \gHom_A(-,A)\circ(-\otimes_A\Omega[d-1])\circ \gHom_{A^o}(-,A)&\cong \gHom_A(\Omega,A)[1-d]\otimes_A-\\
       &\cong \Omega[d-1]\otimes_A-
    \end{align*}
    as functors over $\underline{\MCM}(A^o)$.
    It follows from \cite[Theorem 1.3]{Ue1} that $A^o$ is a noncommutative isolated singularity.
\end{proof}

So, if $A$ is a commonly graded AS-Gorenstein isolated singularity, it is a noncommutative isolated singularity on both sides.

\begin{lemma}\label{global dimension of tails A}
    If $A$ is a commonly graded AS-Gorenstein isolated singularity of dimension $d\geqslant 2$, then $\gldim(\tails A)=d-1$.
\end{lemma}
\begin{proof}
    It follows from a similar argument as in the connected graded case \cite[Corollay 4.5]{Ue1}.
\end{proof}

The following lemma shows a crucial property of MCM modules over commonly graded AS-Gorenstein isolated singularities. It is proved in \cite[Lemma 5.7]{Ue1} for connected graded algebras. A different proof is given in the following.

\begin{lemma}\label{gExt(M,N) is f.d when M is MCM}
Let $A$ be a commonly graded AS-Gorenstein isolated singularity. Then $\gExt_A^i(M,N)$ is finite-dimensional for any $M\in \MCM(A)$, $N\in \gr A$ and $i>0$.
\end{lemma}
\begin{proof}
Suppose the global dimension of $\qgr A$ is $n$. By Proposition \ref{MCM and syzygy},
    \[\gExt_A^i(M,N)\cong \gExt_A^{i+n}(M,\Omega^{n}N).\]
By Proposition \ref{facts about chi condition}, $\gExt_A^{i+n}(M,\Omega^nN)\to \gExt_{\cA}^{i+n}(\cM,\pi\Omega^nN)$ has bounded-above kernel and cokernel. Since $\gExt_{\cA}^{i+n}(\cM,\pi\Omega^nN)=0$, $\gExt_A^{i+n}(M,\Omega^nN)$ is bounded-above. It follows that $\gExt_A^{i+n}(M,\Omega^nN)$ and $\gExt_A^i(M,N)$ are finite-dimensional.
\end{proof}

\subsection{MCM Generators and modulo-torsion-invertible bimodules}
In this subsection, MCM generators over noetherian commonly graded AS-Gorenstein algebras are studied. An MCM $A$-module $M$ is called an \textit{MCM generator} over $A$ if $M$ is  a generator of $\gr A$. In fact, $M$ is a generator of $\gr A$ if and only if it is a generator of $\MCM(A)$ when $A$ is an MCM $A$-module.

Every noetherian commonly graded AS-Gorenstein algebra $A$ of dimension $d\geqslant 2$ is a noncommutative projective coordinate ring. Next proposition shows that under some assumptions the image of every MCM generator in $\qgr A$ is a structure sheaf of $(\qgr A,s)$. Hence MCM generators are modulo-torsion-invertible. Its proof follows from \cite[Theorem 3.10]{Ue2} and \cite[Theorem 2.5]{MU} where the assumptions are slightly stronger than here. Proposition \ref{MCM generator is ample} is also a generalization of \cite[Proposition 3.4]{HY}.

\begin{proposition}\label{MCM generator is ample}
    Suppose $A$ is a noetherian commonly graded AS-Gorenstein algebra of dimension $d\geqslant 2$, $M$ is an MCM generator over $A$. If $\gExt_A^1(M,N)$ is finite-dimensional for any $N\in \MCM(A)$, then $B=\gEnd_A(M)$ is right noetherian and $\gHom_A(M,-)$ induces an equivalence: $\pi_B\gHom_{\cA}(\cM,-):(\qgr A,s)\to (\qgr B,s)$.
\end{proposition}
\begin{proof}
    Since $M$ is an MCM $A$-module, $\depth_A(M)=d\geqslant 2$. We claim that $(\cM,s)$ is ample in $\qgr A$.

    For any $\cN\in \qgr A$, $b_l(NJ^n)$ can be large enough when $n\gg 0$ by Lemma \ref{J^nM and M>n cofinal}, and $\pi (NJ^n)\cong \cN$.
    Since $M_A$ is a generator, there are positive integers $r_1,\cdots,r_p$ such that $\oplus_{i=1}^pM(-r_i)\to NJ^n$ is surjective. So, $\oplus_{i=1}^p\cM(-r_i)\to \cN$ is an epimorphism by the exactness of $\pi$.

    Suppose $f:\cN_1\to \cN_2$ is an epimorphism in $\qgr A$. Consider the  exact sequence $0\to\cK\to\cN_1\to\cN_2\to 0$ in $\qgr A$.

    By \cite[Theorem 2.5]{Ue2}, there exists an exact sequence $0\to L\to Z\to K\to0$ in $\gr A$ where $Z$ is in $\MCM (A)$ and the injective dimension $\idim L<\infty$ such that $\gExt_A^1(M,L)=\gExt_A^2(M,L)=0$.
   Then $\gExt_A^1(M,Z) \cong \gExt_A^1(M,K)$.
    So $\gExt_A^1(M,K)$ is finite-dimensional by assumption. Hence $\gExt_{\cA}^1(\cM,\cK)$ is bounded-above by Proposition \ref{facts about chi condition}.
    By applying $\gHom_{\cA}(\cM,-)$ to the exact sequence $0\to\cK\to\cN_1\to\cN_2\to 0$, it follows that
    \[\gHom_{\cA}(\cM,\cN_1)\to \gHom_{\cA}(\cM,\cN_2)\to \gExt_{\cA}^1(\cM,\cK)\]
    is exact.
    Since $\gExt_{\cA}^1(\cM,\cK)$ is bounded-above, 
    $\gHom_{\cA}(\cM,\cN_1)_{\geqslant n}\to \gHom_{\cA}(\cM,\cN_2)_{\geqslant n}$ is surjective for some integer $n$. So $(\cM,s)$ is ample for $\qgr A$.

    By Proposition \ref{M to (B,F)}, $B=\gEnd_A(M)$ is a noncommutative projective coordinate ring. So $B$ is right noetherian and $\gHom_A(M,-)$ induces an equivalence: $\pi_B\gHom_{\cA}(\cM,-):(\qgr A,s)\to (\qgr B,s)$.
\end{proof}

\begin{corollary}\label{MCM generator is ample for isolated singularity}
Suppose $A$ is a noetherian commonly graded AS-Gorenstein isolated singularity of dimension $d\geqslant 2$. For any MCM generator $M$ over $A$, $B=\gEnd_A(M)$ is noetherian and $M$ is a modulo-torsion-invertible $(B,A)$-bimodule.

    In particular, for any MCM $A$-module $M$, $M\oplus A$ is always an MCM generator over $A$, so $B=\gEnd_A(M\oplus A)$ is noetherian and $M\oplus A$ is a modulo-torsion-invertible $(B,A)$-bimodule.
\end{corollary}
\begin{proof}
     Let $M'=\gHom_A(M,A)$. By Proposition \ref{Hom(-,A) induce duality on MCM}, $M'$ is an MCM $A^o$-module and $B^o\cong \gEnd_{A^o}(M')$. Since $M_A$ is a generator, ${}_AM'$ is also a generator. So $M'$ is an MCM generator over $A^o$.

     Note that the condition in Proposition \ref{MCM generator is ample} holds for $M$ and $M'$ by Lemma \ref{gExt(M,N) is f.d when M is MCM} and $A$ is MCM over itself on both sides. The conclusion follows from Proposition \ref{MCM generator is ample} and Theorem \ref{morita theory}.
\end{proof}

This corollary offers a lot of examples of noncommutative quasi-projective spaces admitting different structure sheaves.

It is proved in \cite[Theorem 3.2]{HY} that $\gEnd_A(M)$ is right noetherian for any MCM $A$-module $M$ if $A$ is a right noetherian commonly graded AS-Gorenstein algebra and $A$ is also a noncommutative isolated singularity.

MCM generators possess nice properties as given in the following proposition, which is proved in \cite{Ue2} for MCM modules containing $A$ as a direct summand where $A$ is a connected graded algebra.

\begin{proposition}\label{property of MCM generator}
    Suppose that $A$ is a noetherian commonly graded AS-Gorenstein algebra, $M$ an MCM generator over $A$. Let $B=\gEnd_A(M)$. If  $\gExt_A^1(M,X)$ and $\gExt_A^1(X,M)$ are finite-dimensional for any $X\in \MCM(A)$, then $B$ is noetherian, and for any $N\in \gr A$,
    \begin{enumerate}
        \item [(1)] $\gHom_A(M,N)$ is a finitely generated $B$-module;
        \item [(2)] there are integers $s_1,\cdots,s_t$ and surjective morphism $f:\oplus_{i=1}^tM(s_i)\to N$ such that
        \[f_*: \, \gHom_A(M,\oplus_{i=1}^tM(s_i)) \to \gHom_A(M,N)\]
        is surjective;
        \item [(3)] the natural morphism $\gHom_A(M,N)\otimes_BM\to N$ is an isomorphism.
    \end{enumerate}
\end{proposition}
\begin{proof}
    By Proposition \ref{MCM generator is ample}, $B$ is right noetherian. Dually, by Proposition \ref{Hom(-,A) induce duality on MCM} and a left module version of Proposition \ref{MCM generator is ample}, $B$ is left noetherian.

    (1) By \cite[Theorem 2.5]{Ue2}, there exists an exact sequence $0\to L\to Z\to N\to0$ in $\gr A$ where $Z$ is in $\MCM (A)$ and the injective dimension $\idim L<\infty$ such that $\gExt_A^1(M,L)=0$. So there is an exact sequence of graded $B$-modules
    \[0\to\gHom_A(M,L)\to\gHom_A(M,Z)\to \gHom_A(M,N)\to 0.\]

    Since $Z$ is an MCM $A$-module, by Proposition \ref{Hom(-,A) induce duality on MCM} there is an injective morphism $Z\to P$ where $P$ is a finitely generated free $A$-module. Since $M$ is a generator in $\gr A$, there is an injective morphism $P\to \tilde{M}$ where $\tilde{M}$ is a finite direct sum of shifts of $M$. Therefore we have an injective morphism of graded $B$-modules
    \[\gHom_A(M,Z)\to\gHom_A(M,\tilde{M}).\]
    Since $B$ is right noetherian 
    and $\gHom_A(M,\tilde{M})$ is a finite direct sum of graded free $B$-modules, $\gHom_A(M,Z)$ is finitely generated. Consequently, $\gHom_A(M,N)$ is a finitely generated graded $B$-module.

    (2) Let $\{f_1,\cdots,f_t\}$ be a set of homogeneous generators of $\gHom_A(M,N)$ as graded $B$-module.
    Set $s_i=\deg f_i$ and $f=(f_i):\oplus_{i=1}^tM(s_i)\to N$. Then for any $g\in\gHom_A(M,N)$, $g$ factors through $f$, that is, there are morphisms $g_i:M\to M(s_i)$ such that $g=\sum f_ig_i$. So
    \[f_*:\,\gHom_A(M,\oplus_{i=1}^tM(s_i)) \to \gHom_A(M,N)\]
    is surjective. Moreover, since $M$ is a generator in $\gr A$,  $\sum_g\im(g:M\to N)=N$. It follows that $f$ is surjective.

    (3) By (2), there are two exact sequences
    \[0\to K\to \oplus_{i=1}^tM(s_i)\to N\to 0,\]
    \[0\to \gHom_A(M,K)\to \gHom_A(M,\oplus_{i=1}^tM(s_i))\to \gHom_A(M,N)\to 0.\]
    Hence the following diagram is commutative \big(let $\tilde{M}=\oplus_{i=1}^tM(s_i)$\big)
    \begin{scriptsize}
        \[
\xymatrix{
& \gHom_A(M,K)\otimes_B M \ar[r] \ar[d]^{\psi_K} & \gHom_A(M,\tilde{M})\otimes_BM \ar[r] \ar[d]^{\psi} & \gHom_A(M,N)\otimes_B M \ar[r] \ar[d]^{\psi_N} & 0\\
0\ar[r] & K \ar[r] & \tilde{M} \ar[r] & N \ar[r] & 0
}\]
    \end{scriptsize}
where the vertical maps are natural and the rows are exact. Note that $\psi$ is an isomorphism. Since $M_A$ is a generator, $\psi_K$ is surjective. 
It follows that $\psi_N$ is an isomorphism.
\end{proof}

Graded endomorphism rings of MCM generators over a commonly graded AS Gorenstein isolated singularity have nice properties as given in the following proposition.

\begin{proposition}\label{endomorphism ring of MCM generator}
    Suppose $A$ is a commonly graded AS-Gorenstein isolated singularity of dimension $d\geqslant 2$ and $M$ is an MCM generator over $A$. Let $B=\gEnd_A(M)$.
    \begin{enumerate}
        \item [(1)] $B$ is noetherian and $B$ admits a balanced dualizing complex.
        \item [(2)] Both the cohomological dimensions of $\Gamma_B$ and $\Gamma_{B^o}$ are $d$.
        \item [(3)] $\gHom_A(M,-)$ induces an equivalence between a full subcategory of $\MCM (A)$ and $\MCM(B)$:
        \[\{N\in\MCM(A)\mid \gExt_A^i(M,N)=0,\forall 0<i<d-1\}\to \MCM(B)\]
         with the quasi-inverse $\gHom_B(M',-)$ where $M'=\gHom_A(M,A)$.
    \end{enumerate}
\end{proposition}
\begin{proof}
    (1) By Corollary \ref{MCM generator is ample for isolated singularity}, $B$ is noetherian and $M$ is a modulo-torsion-invertible $(B,A)$-bimodule.
    Hence $(\cM,s)$ is ample in $\qgr A$ by Theorem \ref{modulo-torsion-invertible bimodule and equivalent functor}.
   It follows from  Theorem \ref{noncommutative Serre theorem for commonly graded algebra} that
    \[\pi_B\gHom_{\cA}(\cM,-):\qgr A\rightleftarrows \qgr B:-\otimes_{\cB}\cM\]
    is an equivalence.

    For any $Y\in \gr B$ and $i\geqslant 2$, by Proposition \ref{cohomology in tails},
    \[R^i\Gamma_B(Y)=\gExt_{\cB}^{i-1}(\cB,\cY)\cong \gExt_{\cA}^{i-1}(\cM,\pi_A(Y\otimes_BM)).\]
    Since $\gExt_A^{i-1}(M,Y\otimes_BM)$ is finite-dimensional by Lemma \ref{gExt(M,N) is f.d when M is MCM}, and the natural morphism $\gExt_A^{i-1}(M,Y\otimes_BM)\to \gExt_{\cA}^{i-1}(\cM,\pi_A(Y\otimes_BM))$ has bounded-above kernel and cokernel, $\gExt_{\cA}^{i-1}(\cM,\pi_A(Y\otimes_BM))$ is bounded-above. This means for any $i\geqslant 2$, $R^i\Gamma_B(Y)$ is bounded-above. Since $B$ satisfies $\chi_1$, $B$ satisfies $\chi$ by Proposition \ref{facts about chi condition}.

    Since $\gldim(\qgr A)=d-1$,
    $R^{d+1}\Gamma_B(Y)\cong\gExt_{\cA}^d(\cM,\pi_A(Y\otimes_BM))=0$ for any $Y\in \gr B$.
    It follows that $R^{d+1}\Gamma_B(Y)=0$ for any $Y\in \Gr B$ as $B$ is right noetherian and $R^i\Gamma_B$ commutes with direct limit. So $\Gamma_B$ has finite cohomological dimension.

    Let $M'=\gHom_A(M,A)$. It follows from the proof of Corollary \ref{MCM generator is ample for isolated singularity} that $M'$ is an MCM generator over $A^o$. So by a dual argument, $B^o$ satisfies $\chi$ and $\Gamma_{B^o}$ has finite cohomological dimension. By a commonly graded version of \cite[Theorem 6.3]{V3}, $B$ admits a balanced dualizing complex.

    (2) Note that
    $R^d\Gamma_B(M')\cong\gExt_{\cB}^{d-1}(\cB,\cM')\cong \gExt_{\cA}^{d-1}(\cM,\cA).$

    Let $0\to A_A\to I^0\xrightarrow[]{\partial^0} I^1\xrightarrow[]{\partial^1}\cdots$ be the minimal graded injective resolution of $A_A$. By Theorem \ref{minimal injective resolution over isolated singularity} and Lemma \ref{global dimension of tails A}, $I^i$ is torsion-free for any $0\leqslant i<d$ and $I^d$ is torsion. By \cite[7.1.5]{AZ}, $\gExt_{\cA}^{d-1}(\cM,\cA)$ is isomorphic to the $(d-1)$-th cohomology of the complex
    \[0\to\gHom_A(M,I^0)\xrightarrow[]{\partial^0_*}  \cdots\to\gHom_A(M,I^{d-2})\xrightarrow[]{\partial^{d-2}_*}\gHom_A(M,I^{d-1})\to 0.\]

    Since $M\in\MCM(A)$, $\gExt_A^{d-1}(M,A)=\gExt_A^d(M,A)=0$. So
    \[\gHom_A(M,I^{d-2})\xrightarrow[]{\partial^{d-2}_*}\gHom_A(M,I^{d-1})\xrightarrow[]{\partial^{d-1}_*}\gHom_A(M,I^d)\to 0\]
    is exact. Since $M_A$ is a generator, $\gHom_A(M,I^d)\neq 0$. 
    Hence $\im(\partial^{d-2}_*)=\Ker(\partial^{d-1}_*)\neq \gHom_A(M,I^{d-1})$. Therefore $\gExt_{\cA}^{d-1}(\cM,\cA)\neq 0$ and $R^d\Gamma_B(M')\neq 0$.

    It follows from (1) that the cohomological dimension of $\Gamma_B$ is $d$. Dually, the cohomological dimension of $\Gamma_{B^o}$ is $d$.

   (3) Suppose $N\in \MCM(A)$ and $\gExt_A^i(M,N)=0$ for all $0<i<d-1$. 
   By Theorem \ref{equivalence: modules with depth geqslant 2}, $Y:=\gHom_A(M,N)$ is a finitely generated graded $B$-module and $\depth_B(Y)\geqslant 2$.
   For any $i\geqslant 2$,
   \[R^i\Gamma_B(Y)\cong \gExt_{\cB}^{i-1}(\cB,\cY)\cong \gExt_{\cA}^{i-1}(\cM,\cN).\]
    Since $\depth_A(N)=d$, by Lemma \ref{Ext in tails cong Ext in gr when depth>n}, $\gExt_{\cA}^{i-1}(\cM,\cN)\cong\gExt_A^{i-1}(M,N)=0$ for any $1\leqslant i-1<d-1$. So $R^i\Gamma_B(Y)=0$ for any $i<d$.

    It remains to show $R^d\Gamma_B(Y)\neq 0$, which in fact follows from a similar argument as (2) by replacing $M'$ with $Y$. So $\depth_B(Y)=d$ and $\gHom_A(M,N)$ is an MCM $B$-module.

    Conversely, suppose $Y\in \MCM(B)$. By Theorem \ref{equivalence: modules with depth geqslant 2}, there is a finitely generated graded $A$-module $N$ with $\depth_A(N)\geqslant 2$ such that $Y\cong \gHom_A(M,N)$.

    Since $M_A$ is a generator, $M'_B$ is a finitely generated projective $B$-module. By Proposition \ref{cohomology in tails} and Lemma \ref{Ext in tails cong Ext in gr when depth>n}, for any $2\leqslant i<d$,
    \[R^i\Gamma_A(N)\cong \gExt_{\cA}^{i-1}(\cA,\cN)\cong \gExt^{i-1}_{\cB}(\cM',\cY)\cong \gExt_B^{i-1}(M',Y)=0.\]

    Since $A$ is a noncommutative isolated singularity, $\depth_A(N)=d$ by Theorem \ref{minimal injective resolution over isolated singularity}. So $N$ is an MCM $A$-module, and by Lemma \ref{Ext in tails cong Ext in gr when depth>n}, for any $0<i<d-1$,
    \[\gExt_A^i(M,N)\cong \gExt_{\cA}^i(\cM,\cN)\cong \gExt^i_{\cB}(\cB,\cY)\cong R^{i+1}\Gamma_B(Y)=0.\]

    The assertion follows from Theorem \ref{equivalence: modules with depth geqslant 2}.
\end{proof}

\begin{proposition}\label{when B is balanced CM}
     Suppose $A$ is a commonly graded AS-Gorenstein isolated singularity of dimension $d\geqslant 2$ and $B$ is a noncommutative projective coordinate ring. If $M$ is a modulo-torsion-invertible $(B,A)$-bimodule with $\depth_AM\geqslant 2$ and $M$ is a generator in $\gr A$, then the following are equivalent.
     \begin{enumerate}
         \item [(1)] $M$ is an MCM $A$-module and $\gExt_A^i(M,M)=0$ for all $0<i<d-1$.
         \item [(2)] $B$ is a balanced Cohen-Macaulay algebra of dimension $d$.
     \end{enumerate}
\end{proposition}
\begin{proof}
    Suppose $M$ is an MCM $A$-module and $\gExt_A^i(M,M)=0$ for all $0<i<d-1$. By Proposition \ref{endomorphism ring of MCM generator} (1), $B$ is a noetherian algebra and admits a balanced dualizing complex. By Proposition \ref{endomorphism ring of MCM generator} (3), $B$ is an MCM $B$-module. Therefore $B$ is a balanced Cohen-Macaulay algebra of dimension $d$.

    Conversely, suppose $B$ is a balanced Cohen-Macaulay algebra of dimension $d$. Let $M'=\gHom_A(M,A)$.

    By Proposition \ref{cohomology in tails}, for any $2\leqslant i<d$,
    \[R^i\Gamma_A(M)\cong\gExt_{\cA}^{i-1}(\cA,\cM)\cong \gExt_{\cB}^{i-1}(\cM',\cB).\]

    Since $M_A$ is a generator, $M'_B$ is a finitely generated projective $B$-module.
    Thus by Lemma \ref{Ext in tails cong Ext in gr when depth>n}, for any $2\leqslant i<d$,
    \[\gExt_{\cB}^{i-1}(\cM',\cB)\cong \gExt_B^{i-1}(M',B)=0.\]
    Therefore $R^i\Gamma_A(M)=0$. Moreover, it follows from $\depth_A(M)\geqslant 2$ that $M$ is an MCM $A$-module.
    By Proposition \ref{endomorphism ring of MCM generator}, $\gExt_A^i(M,M)=0$ for any $0<i<d-1$.
\end{proof}

\section{Noncommutative resolutions}\label{Noncommutative Crepant Resolutions}

In this section, we define and study noncommutative resolutions for commonly graded AS-Gorenstein isolated singularities. The noncommutative resolutions in Definition \ref{definition of NCCR} are always commonly graded Artin-schelter regular algebras (see Theorem \ref{NCCR is AS regular}), and are strongly related to cluster-tilting modules (see Theorem \ref{d-1 CT module and NCCR}).
A version of the noncommutative Bondal-Orlov conjecture is proved to be true in dimensions $2$ and $3$ (see Theorem \ref{proof for BO conjecture}).

\subsection{Noncommutative  Resolutions}
The following definition is motivated by the one in \cite{V1,IR,IW1} etc.

\begin{definition}\label{definition of NCCR}
    Let $A$ be a commonly graded AS-Gorenstein isolated singularity of dimension $d\geqslant 2$.  A \textit{noncommutative resolution} of $A$ is a graded ring $B=\gEnd_A(M)$ for some MCM generator $M_A$ such that $B$ has finite graded global dimension $d$ and $B_B$ is an MCM $B$-module. In this case, the  noncommutative resolution $B$ is said to be \textit{given by $M_A$}.
\end{definition}

Next proposition follows from Corollary \ref{MCM generator is ample for isolated singularity}.
\begin{proposition}\label{tails of NCCR equi. to tails of singularity}
     Suppose $A$ is a commonly graded AS-Gorenstein isolated singularity of dimension $d\geqslant 2$ and $M$ is an MCM generator over $A$. If $M_A$ gives a noncommutative resolution $B=\gEnd_A(M)$ of $A$, then the graded Morita context defined by $M_A$ induces an equivalence
\[-\otimes_{\cA}\cM':(\qgr A,s)\rightleftarrows (\qgr B,s):-\otimes_{\cB}\cM\]
where $M'=\gHom_A(M,A)$.
\end{proposition}

In fact, any noncommutative resolution $B$ of $A$ is a noetherian commonly graded AS-regular algebra of dimension $d$, which will be proved in Theorem \ref{NCCR is AS regular}.

As mentioned before, AS-regular algebras are regarded as the coordinate rings of noncommutative projective spaces, and noncommutative projective schemes associated to AS-regular algebras are smooth. To give a  noncommutative resolution of a noncommutative isolated singularity $A$ is to find a commonly graded AS-regular algebra $B$ such that the noncommutative quasi-projective space $(\qgr A,s)$ is equivalent to the noncommutative projective space $(\qgr B,s)$. It follows that if $A$ admits a noncommutative resolution then $(\qgr A,s)$ is a noncommutative projective space (see Definition \ref{noncommutative-proj-space}). Proposition \ref{tails of NCCR equi. to tails of singularity} and Theorem \ref{NCCR is AS regular} justify this.

\begin{remark}\label{explain the name NCCR}
In algebraic geometry, a crepant resolution $f:Y\to X$ of $X$ 
is a proper birational morphism $f$ such that $Y$ is smooth and $f^*\omega_X=\omega_Y$, where $\omega_X$ and $\omega_Y$ are the corresponding canonical bundles.

If $M$ in Definition \ref{definition of NCCR} is an $(A,A)$-bimodule, then there is a natural map $\varphi:A\to \gEnd_A(M)=B$. If ${}_AM_A$ is finitely generated on both sides, then ${}_AB$ and $B_A$ are finitely generated. This holds automatically if $A$ is commutative.  In noncommutative realm, ``proper" is replaced by that both ${}_AB$ and $B_A$ are finitely generated via $\varphi$. ``Nonsingularity" is replaced by $\gldim B=d$. Noncommutative resolutions of the invariant subrings of Hopf actions are such examples, where the noncommutative resolution of $A$ is given by an $(A,A)$-bimodule finitely generated on both sides (see Example \ref{example of Zhu and CYZ}).

As explained in \cite{V1}, ``birationality" should be replaced by Morita equivalence in noncommutative geometry. In noncommutative projective geometry, ``birationality" should be replaced by the equivalence of noncommutative quasi-projective spaces. By the Morita theory of the noncommutative quasi-projective spaces established in \S \ref{Morita Theory for Noncommutative quasi-projective spaces}, $B$ is of the form $\gEnd_A(M)$ where $M$ is a modulo-torsion-invertible $(B,A)$-bimodule with $\depth_AM\geqslant 2$.

Let $\Omega_A$ be the balanced dualizing module of $A$. It follows from Proposition \ref{local cohomology over B and over A when B_A finite} that the balanced dualizing complex of $B$ satisfies $D(R\Gamma_B(B))\cong D(R\Gamma_A(B))\cong R\gHom_A(B,\Omega_A)$. The condition that the balanced dualizing complex of $B$ is concentrated in degree $d$ is used to substitute for ``crepancy", which is equivalent to say that $B_A$ is MCM. This is the case if the modulo-torsion-invertible $(B,A)$-bimodule $M$ is an MCM generator over $A$ and $\gExt_A^i(M,M)=0$ for all $0<i<d-1$ by Proposition \ref{when B is balanced CM}.
\end{remark}

Next we prove that any noncommutative resolution in the sense of Definition \ref{definition of NCCR} is a commonly graded AS-regular algebra. The following lemma is well-known in connected graded case.

\begin{lemma}\label{balanced CM algebra with finite gldim}
    Let $B$ be a balanced Cohen-Macaulay algebra of dimension $d$. If $\gldim B=d$, then $B$ is a commonly graded AS-regular algebra.
\end{lemma}
\begin{proof}
   It follows from $\gldim B=d$ that for any graded simple $B$-module $Y$ the following spectral sequence converges
    \[E^{p,q}_2=\gExt_{B^o}^p(\gExt_B^{-q}(Y,B),B)\Rightarrow\left\{
    \begin{aligned}
        0,\quad & p+q\neq 0,\\
        Y,\quad & p+q=0.
    \end{aligned}
    \right.
    \]
    Since $\depth_B(B)=\depth_{B^o}(B)=d$,
    $\gExt_{B^o}^d(\gExt_B^d(Y,B),B)\cong Y.$

    Similarly, for any graded simple $B^o$-module $\tilde{Y}$,
    $\gExt_B^d(\gExt_{B^o}^d(\tilde{Y},B),B)\cong \tilde{Y}.$ Thus $B$ satisfies the condition $(C)$ and $(C^o)$ in \cite[Theorem 5.15]{LW}. Hence $B$ is a commonly graded AS-regular algebra.
\end{proof}

\begin{theorem}\label{NCCR is AS regular}
    Suppose $A$ is a commonly graded AS-Gorenstein isolated singularity of dimension $d\geqslant 2$. Then the following are equivalent.
    \begin{enumerate}
        \item [(1)] $B$ is a noncommutative resolution of $A$ given by $M_A$.
        \item [(2)] $B^o$ is a noncommutative resolution of $A^o$ given by ${}_AM'=\gHom_A(M,A)$ and $M_A$ is an MCM module.
        \item [(3)] $B=\gEnd_A(M)$ is a noetherian commonly graded AS-regular algebra of dimension $d$ and $M_A$ is an MCM generator.
    \end{enumerate}
\end{theorem}
\begin{proof}
(1) $\Rightarrow$ (3)
    By Proposition \ref{endomorphism ring of MCM generator}, $B$ is noetherian and $B$ admits a balanced dualizing complex. Since $B$ is an MCM $B$-module, $B$ is balanced Cohen-Macaulay of dimension $d$. It follows from Lemma \ref{balanced CM algebra with finite gldim} that $B$ is a commonly graded AS-regular algebra of dimension $d$.

(3) $\Rightarrow$ (1) It follows from the definition.

(2) $\Leftrightarrow$ (3) Note that ${}_AM'$ is an MCM generator. By Proposition \ref{Hom(-,A) induce duality on MCM}, $B\cong\gEnd_A(M)\cong\gEnd_{A^o}(M')$. So the proof is same as (1) $\Leftrightarrow$ (3).
\end{proof}

\subsection{Noncommutative resolutions versus cluster tilting modules}
Similar to the case of module-finite algebras \cite{I2}, the following theorem reveals that every $(d-1)$-cluster tilting module gives a noncommutative resolution and every noncommutative resolution is given by a $(d-1)$-cluster tilting module. Moreover, for an MCM generator $M_A$ over a commonly graded AS-Gorenstein isolated singularity $A$, part of the conditions in the definition of $(d-1)$-cluster tilting modules is enough to assure $M_A$ being $(d-1)$-cluster tilting. It follows from Definition \ref{cluster tilting module} that
every cluster tilting module is an MCM generator.

\begin{theorem}\label{d-1 CT module and NCCR}
Suppose $A$ is a commonly graded AS-Gorenstein isolated singularity of dimension $d\geqslant 2$.  Then the following are equivalent.
\begin{enumerate}
    \item [(1)] $M$ is a $(d-1)$-cluster tilting $A$-module.
    \item [(2)] $\add_A M=\{N\in \MCM (A)\mid \gExt_A^i(M,N)=0,\forall\, 0<i<d-1\}$.
    \item [(3)] $\add_A M=\{N\in \MCM (A)\mid \gExt_A^i(N,M)=0,\forall\, 0<i<d-1\}$.
    \item [(4)] $B=\gEnd_A(M)$ is a noncommutative resolution of $A$ given by $M_A$.
\end{enumerate}
\end{theorem}
\begin{proof}
    Note that all (1), (2), (3), (4) imply that $M_A$ is an MCM generator.

    (1) $\Leftrightarrow$ (2) and (3) By Definition \ref{cluster tilting module}.

    (2) $\Rightarrow$ (4) It follows from Proposition \ref{when B is balanced CM} that $B$ is a balanced Cohen-Macaulay algebra of dimension $d$.

    By the projectivisation (say see \cite[\uppercase\expandafter{\romannumeral6} Lemma 3.1]{SS}), there is an equivalence
    \[\gHom_A(M,-):\add_A(M)\to \proj B\]
    where $\proj B$ is the category of finitely generated graded projective $B$-modules.
    It follows from Proposition \ref{endomorphism ring of MCM generator} that $\MCM(B)=\proj B$.

    For any $Y\in\gr B$, let
    \[0\to K\to P_{d-1}\xrightarrow[]{\partial_{d-1}}\cdots\to P_1\xrightarrow[]{\partial_1} P_0\xrightarrow[]{\partial_0} Y\to 0\]
    be a finitely generated graded resolution of $Y$ where $P_i$ is graded projective and $K=\Ker\partial_{d-1}$. By Lemma \ref{depth in exact sequences}, $\depth_B(K)\geqslant d$.
    Since the cohomological dimension of $\Gamma_B$ is $d$, $\depth_B(K)=d$ by Lemma \ref{depth < infinity}. Hence $K$ is an MCM $B$-module. Since $\MCM(B)=\proj B$, $K$ is a finitely generated graded projective $B$-module. It follows that the projective dimension of $Y$ is no more than $d$. Since the cohomological dimension of $\Gamma_B$ is $d$, $\gldim B=d$.

    It follows from Lemma \ref{balanced CM algebra with finite gldim} that $B$ is a noncommutative resolution of $A$ given by $M$.

    (4) $\Rightarrow$ (2) By Theorem \ref{NCCR is AS regular}, $B$ is a commonly graded AS-regular algebra of dimension $d$. Hence by \cite[Theorem 1.7]{HY}, $\MCM(B)=\proj B$.

    Note that $\add_AM$ and $\{N\in\MCM(A)\mid \gExt_A^i(M,N)=0,\forall \, 0<i<d-1\}$ are contained in $\{N\in \gr A\mid \depth_A(N)\geqslant 2\}$. By Theorem \ref{equivalence: modules with depth geqslant 2}, $\gHom_A(M,-)$, when restricted to $\{N\in \gr A\mid \depth_A(N)\geqslant 2\}$, gives an equivalent functor.
    By the projectivisation and Proposition \ref{endomorphism ring of MCM generator}, $\gHom_A(M,-)$ also gives the following equivalences
    \[\add_AM\to \proj B=\MCM(B)\]
    and
    \[\{N\in\MCM(A)\mid \gExt_A^i(M,N)=0,\forall\, 0<i<d-1\}\to \MCM(B),\]
    so
    \[\add_A M=\{N\in\MCM(A)\mid \gExt_A^i(M,N)=0,\forall\, 0<i<d-1\}.\]

    (4) $\Rightarrow$ (3) By Theorem \ref{NCCR is AS regular}, $B^o$ is a noncommutative resolution of $A^o$ given by ${}_AM'=\gHom_A(M,A)$. By ``(4) $\Rightarrow$ (2)",
    \[\add_{A^o}M'=\{\tilde{N}\in \MCM (A^o)\mid \gExt^i_{A^o}(M',\tilde{N})=0,\forall\, 0<i<d-1\}.\]
    Since $\gHom_{A^o}(-,A)$ gives a duality between $\MCM(A^o)$ and $\MCM(A)$, it follows from $\gExt_A^j(X,A)=0$ for any $X\in \MCM A$ and $j>0$ that
    \[\add_A M=\{N\in \MCM (A)\mid \gExt_A^i(N,M)=0,\forall\, 0<i<d-1\}.
    \qedhere\]
\end{proof}

\begin{corollary}\label{M is CT iff M' is CT}
If $M_A$ is an MCM module, then $M$ is a $(d-1)$-cluster tilting $A$-module if and only if $\gHom_A(M,A)$ is a $(d-1)$-cluster tilting $A^o$-module.
\end{corollary}

\subsection{Uniqueness of the noncommutative resolutions}

In this final subsection, we prove that a version of the noncommutative Bondal-Orlov conjecture \cite{V1} is true in dimension $2$ and dimension $3$ cases. The key of the proof is that the endomorphism rings of $1$-cluster tilting modules are Morita equivalent, and the endomorphism rings of $2$-cluster tilting modules are derived Morita equivalent as in \cite[Theorem 5.3.2]{I2}. 
When the dimension is greater than $3$, this noncommutative version of the Bondal-Orlov conjecture is also open.  But by Corollary \ref{center of noncommutative coordinate ring}, the centers of all noncommutative resolutions of a commonly graded AS-Gorenstein isolated singularity $A$ are isomorphic, and are isomorphic to the center of $A$.

\begin{theorem}\label{proof for BO conjecture}
Suppose $A$ is a commonly graded AS-Gorenstein isolated singularity of dimension $d$.
\begin{enumerate}
    \item [(1)] If $d=2$, then all  noncommutative resolutions of $A$ are Morita equivalent.
    \item [(2)] If $d=3$, then all  noncommutative resolutions of $A$ are derived Morita equivalent.
\end{enumerate}
\end{theorem}
\begin{proof}
    (1) This is actually proved in \cite[Theorem 2.3]{HY}.

    (2) Suppose $B_1$ and $B_2$ are noncommutative resolutions of $A$ given by $M_1$ and $M_2$ respectively. By Theorem \ref{d-1 CT module and NCCR}, $M_1$ and $M_2$ are $2$-cluster tilting $A$-modules.
    Let ${}_{B_2}U_{B_1}=\gHom_A(M_1,M_2)$. We claim that $U$ is a tilting $(B_2,B_1)$-bimodule, which induces the derived Morita equivalence between $B_1$ and $B_2$.
	
It follows from Proposition \ref{equivalence: modules with depth geqslant 2} that
\begin{align*}
\gEnd_{B_1}(U)&=\gHom_{B_1}(\gHom_A(M_1,M_2),\gHom_A(M_1,M_2))\\
&\cong \gHom_A(M_2,M_2)\\
&\cong B_2.
\end{align*}
	
By Proposition \ref{property of MCM generator}, there is an exact sequence
	\[0\rightarrow M\xrightarrow{g} \oplus M_1(s_i)\xrightarrow[]{f} M_2\rightarrow 0\]
such that
\[f_*:\gHom_A(M_1,\oplus M_1(s_i))\rightarrow \gHom_A(M_1,M_2)\]
is surjective. Since both $M_1$ and $M_2$ are MCM $A$-modules, so is $M$.

Since $\gExt_A^1(M_1,M_1)=0$, it follows from the long exact sequence that
 \[\gExt_A^1(M_1,M)=0.\]
 Since $M_1$ is a $2$-cluster tilting $A$-module, $M\in \add_A M_1$. Thus
 \begin{equation}\label{prroj-reso-of-U}
     0\to \gHom_A(M_1,M)\to \gHom_A(M_1,\oplus M_1(s_i))\rightarrow \gHom_A(M_1,M_2)\to 0 \tag{$\star\star$}
 \end{equation}
 is a graded projective resolution of $U_{B_1}$ and  $\pdim_{B_1} U\leqslant 1$.

 By using the graded projective resolution \eqref{prroj-reso-of-U} and  Proposition \ref{equivalence: modules with depth geqslant 2}, $\gExt_{B_1}^j(U,U)$ is isomorphic to the $j$-th cohomology of the following complex:
 \[0\to X^0=\gHom_A(\oplus M_1(s_i),M_2)\xrightarrow[]{g^*} X^1=\gHom_A(M,M_2) \xrightarrow[] 0 \cdots\]
Since $\gExt_A^1(M_2,M_2)=0$, $g^*$ is surjective. So $\gExt_{B_1}^1(U,U)=0$.

Let $M_1'=\gHom_A(M_1,A)$ and $M_2'=\gHom_A(M_2,A)$. Then $M_2'$ is a $2$-cluster tilting $A^o$-module by Corollary \ref{M is CT iff M' is CT}.

By Proposition \ref{property of MCM generator}, there is an exact sequence
\[0\to \tilde{M}\to \oplus M_2'(\tilde{s}_i)\to M_1'\to 0.\]
Similar to previous proof, $\tilde{M}\in\add_{A^o}M_2'$. By Proposition \ref{Hom(-,A) induce duality on MCM}, the following sequence
\[0\to M_1 \to \oplus M_2(-\tilde{s}_i)\to \gHom_{A^o}(\tilde{M},A)\to 0\]
is exact and $\gHom_{A^o}(\tilde{M},A)$ is in $\add_AM_2$. Let $N=\gHom_{A^o}(\tilde{M},A)$.
Since $\gExt_A^1(M_1,M_1)=0$,
\[0\to B_1=\gHom_A(M_1,M_1)\to \gHom_A(M_1,\oplus M_2(-\tilde{s}_i))\to \gHom_A(M_1,N)\to 0\]
is exact and $\gHom_A(M_1,\oplus M_2(-\tilde{s}_i))$, $ \gHom_A(M_1,N)$ are in $\add_{B_1}U$.

In conclusion, $U$ is a tilting $(B_2,B_1)$-bimodule and $B_1$ and $B_2$ are derived equivalent by \cite[Theorem 6.3]{Ri}.
\end{proof}

\begin{example}
Keep the notations as Example \ref{projective space P1}. One sees that $A_{A^G}$ is an $1$-cluster tilting module of $A^G$. By \cite{Au}, there is an isomorphism of algebras $\gEnd_{A^G}(A)\cong A\#G$. So the isomorphism
\[(\qgr A\#G, \pi(A\#G),s)\cong (\qgr A^G,\cA,s)\]
shows $A\#G$ is a noncommutative resolution of $A^G$. Moreover, by Theorem \ref{proof for BO conjecture}, $A\#G$ is the only noncommutative resolution of $A^G$ under Morita equivalence.
\end{example}

Invariant rings of Hopf actions provide a lot of examples of AS-Gorenstein isolated singularities,  the study of their noncommutative resolutions will appear in a subsequent article. Here is a well-studied example in \cite{CYZ,Z}.
\begin{example}\label{example of Zhu and CYZ}
Let $A=k_{-1}[x_0,\cdots,x_{n-1}]$ be the $(-1)$-skew polynomial ring, that is
\[A=k\langle x_0,\cdots,x_{n-1}\rangle/(x_ix_j+x_jx_i,i\neq j).\]
Let $G$ be the cyclic group of order $n$ generated by $\sigma=(0,1,\cdots,n-1)$ that acts on generators by $\sigma\cdot x_i=x_{i+1}$ for all $i\in \mathbb{Z}/n\mathbb{Z}$.
\begin{enumerate}
    \item [(1)] \cite[Theorem 1.3]{Z} If $3 \nmid n$ and $5 \nmid n$, then the $G$-action is ample. Consequently, $A\# G$ is a noncommutative resolution of $A^G$ given by $A_{A^G}$.
    \item [(2)] \cite[Theorem 0.4]{CYZ} If $3 \mid n$  or $5 \mid n$, then the $G$-action is not ample. Consequently, $A\#{G}$ is not a noncommutative resolution of $A^{G}$.
\end{enumerate}
\end{example}

\section*{Acknowledgments}
The authors express their deep appreciation for the invaluable comments and suggestions provided by the referee.


\begin{thebibliography}{AAA}
\bibitem[AS]{AS} M. Artin, W. Schelter, Graded algebras of dimension 3, Adv. Math. 66 (1987), 171-216.

\bibitem[ATV]{ATV} M. Artin, J. Tate, M. Van den Bergh, Some algebras associated to automorphisms of elliptic curves,
The Grothendieck Festschrift, vol. I, Progr. Math., vol. 86, Birkhäuser Boston, Boston, 1990, 33--85.

\bibitem[AV]{AV} M. Artin, M. Van den Bergh, Twisted homogeneous coordinate rings, J. Algebra 133 (1990), no. 2, 249–271.

\bibitem[AZ]{AZ} M. Artin, J.J. Zhang, Noncommutative projective schemes, Adv. Math. 109 (1994), 228-287.

\bibitem[Au]{Au} M. Auslander, On the purity of the branch locus, Amer. J. Math. 84 (1962), 116-125.


\bibitem[BKR]{BKR} T. Bridgeland, A. King, M. Reid, The McKay correspondence as an
equivalence of derived categories, J. Amer. Math. Soc. 14 (2001), 535-554.

\bibitem[BO1]{BO1} A. Bondal, D. Orlov, Derived categories of coherent sheaves, Proceedings of the International Congress of Mathematicians, Vol. II (Beijing, 2002), 47-56, Higher Ed. Press, Beijing, 2002.

\bibitem[BO2]{BO2}  A. Bondal, D. Orlov, Semi-orthogonal decompositions for algebraic varieties, available as alggeom/950601, 1996.

\bibitem[Bri]{Bri} T. Bridgeland, Flops and derived categories, Invent. Math. 147 (2002), 613-632.

\bibitem[CKWZ]{CKWZ2} K. Chan, E. Kirkman, C. Walton, J.J. Zhang, McKay correspondence for semisimple Hopf actions on regular graded algebras, J. Noncommut. Geom. 13 (2019), 87-114.

\bibitem[CYZ]{CYZ} K. Chan, A. Young, J.J. Zhang, Noncommutative cyclic isolated singularities, Trans. Am. Math. Soc. 373 (2020), 4319-4358.

\bibitem[Fa]{Fa} C. Faith, Algebra: Rings, Modules and Categories \uppercase\expandafter{\romannumeral1}, Springer-Verlag, New York, 1973.

\bibitem[Ha]{Ha} R. Hartshorne, Algebraic geometry, Springer-Verlag, New York, 1977.

\bibitem[HY]{HY} J. He, Y. Ye, Pre-resolution of noncommutative isolated singularities, Pacific J. Math. 316 (2022), 367-394.

\bibitem[I1]{I1} O. Iyama, Higher-dimensional Auslander-Reiten theory on maximal orthogonal subcategories, Adv. Math. 210 (2007), 22-50.

\bibitem[I2]{I2} O. Iyama, Auslander correspondence, Adv. Math. 210 (2007) 51-82.

\bibitem[IR]{IR} O. Iyama, I. Reiten, Fomin-Zelevinsky mutation and tilting modules over Calabi-Yau algebras, Amer. J. Math. 130 (2008), 1087-1149.

\bibitem[IW1]{IW1} O. Iyama, M. Wemyss, On the noncommutative Bondal-Orlov Conjecture, J. Reine Angew. Math. 683 (2013), 119-128.

\bibitem[IW2]{IW2} O. Iyama, M. Wemyss, Maximal modifications and Auslander-Reiten duality for nonisolated singularities, Invent. Math. 197 (2014), 521-586.

\bibitem[Ja]{Ja} N, Jacobson, Basic Algebra II,  W. H. Freeman \& Compagny, 1976.

\bibitem[Jo]{Jo}  P. J{\o}rgensen, Finite Cohen-Macaulay type and smooth non-commutative schemes, Canad. J. Math. 60 (2008) 379-390.

\bibitem[KV]{KV} M. Kapranov, E. Vasserot, Kleinian singularities, derived categories
and Hall algebras. Math. Ann. 316 (2000), 565-576.

\bibitem[Leu]{Leu} Graham J. Leuschke, Non-commutative crepant resolutions: scenes From categorical geometry. Progress in commuta-tive algebra, vol. 1, de Gruyter, Berlin, 2012. 293–361.


\bibitem[Lou]{L} Q. Lou, Invariant theory of AS-regular algebras and co-Poisson structures (Chinese), Ph.D. thesis, Fudan University, 2016.

\bibitem[LW1]{LW} H.-N. Li, Q.-S. Wu, Commonly graded algebras and their homological properties, submitted.

\bibitem[LW2]{LW2} H.-N. Li, Q.-S. Wu, Regular $\mathbb{Z}$-graded local rings and graded isolated singularities, Proc. Amer. Math. Soc. 153 (2025), no. 10, 4255–4270.

\bibitem[MU]{MU} I. Mori and K. Ueyama, Ample group actions on AS-regular algebras and noncommutative graded isolated singularities, Trans. Amer. Math. Soc. 368 (2016), 7359-7383.

\bibitem[NO]{NO} C. N$\breve{\mathrm{a}}$st$\breve{\mathrm{a}}$sescu, F. Van Oystaeyen, Graded and Filtered Rings and Modules, Lecture Notes in Mathematics, vol 758. Springer, Berlin, 1979.

\bibitem[Po]{Po} N. Popescu, Abelian categories with applications to rings and modules, L.M.S. Monographs, Academic Press, New York, 1973.

\bibitem[QWZ]{QWZ} X. Qin, Y. Wang, J.J. Zhang, Noncommutative quasi-resolutions, J. Algebra 536 (2019), 102-148.

\bibitem[Ri]{Ri} J. Rickard. Morita theory for derived categories, J. London Math. Soc. (2) 39 (1989), 436-456.

\bibitem[Ro]{Ro} A. L. Rosenberg, Noncommutative schemes, Compositio Math. 112 (1998), 93-125.

\bibitem[RR]{RR} M. L. Reyes, D. Rogalski, Graded twisted Calabi-Yau algebras are generalized Artin-Schelter regular. Nagoya Math. J. 245 (2022), 100-153.

\bibitem[Sm]{S} S. P. Smith, Corrigendum to "Maps between non-commutative spaces'', Trans. Amer. Math. Soc. 368 (2016), 8295-8302.

\bibitem[Se]{Se} J.-P. Serre, Faisceaus alg\'ebriques coh\'erents, Ann. Math. 61 (1955), 197-278.

\bibitem[SS]{SS} D. Simson, A. Skowro\'{e}ski, Elements of the representation theory of associative algebras, London Mathematical Society Student Texts, Cambridge University Press (2007).

\bibitem[Ue1]{Ue1} K. Ueyama, Graded maximal Cohen-Macaulay modules over noncommutative graded Gorenstein isolated singularities, J. Algebra 383 (2013), 85-103.

\bibitem[Ue2]{Ue2} K. Ueyama, Noncommutative graded algebras of finite Cohen-Macaulay representation type, Proc. Amer. Math. Soc. 143 (2015), 3703-3715.

\bibitem[Ue3]{Ue3} K. Ueyama, Cluster tilting modules and noncommutative projective schemes, Pacific J. Math. 289 (2017), 449-468.

\bibitem[Vak]{Vak} The rising sea: foundations of algebraic geometry, 2024, Online notes.

\bibitem[V1]{V1} M. Van den Bergh, Non-commutative crepant resolutions. The legacy of Niels Henrik Abel, 749-770, Springer, Berlin, 2004.

\bibitem[V2]{V2} M. Van den Bergh, Three-dimensional flops and noncommutative rings. Duke Math. J. 122 (2004), 423-455.

\bibitem[V3]{V3} M. Van den Bergh, Existence theorems for dualizing complexes over non-commutative graded and filtered rings, J. Algebra 195 (1997), no.~2, 662-679.

\bibitem[Zhu]{Z} R. Zhu, Some invariant subalgebras are graded isolated singularities, J. Pure Appl. Algebra 226 (2022), Paper No. 107087, 6 pp.
\end{thebibliography}
\end{document}